\documentclass{article}
\usepackage{latexsym,amssymb,amsmath}
\usepackage{graphicx}
\usepackage{color}
\usepackage[emtex,matrix,arrow]{xy}
 \newtheorem{thm}{Theorem}
\newtheorem{cor}[thm]{Corollary}
\newtheorem{ex}[thm]{Example}

\newtheorem{lem}[thm]{Lemma}
\newtheorem{prop}[thm]{Proposition}
\newtheorem{defn}[thm]{Definition}

\numberwithin{equation}{section}

\newcommand{\BB}{\mathbb B}
\newcommand{\Sd}{\mathbb S}

\newcommand{\Integer}{\mathbb Z}
\newcommand{\Natural}{\mathbb N}

\newcommand{\To}{\rightarrow}

\newcommand{\Aa}{\mathcal{A}}

\newcommand{\Gg}{\mathcal{G}}

\newcommand{\Ll}{\mathcal{L}}

\newcommand{\Cc}{\mathcal{C}}

\newcommand{\Dd}{\mathcal{D}}

\newcommand{\F}{\mathbb{F}}
\newcommand{\V}{\mathbb{V}}

\newcommand{\G}{\mathbb{G}}
\newcommand{\Lb}{\mathbb{L}}

\newcommand{\Cgg}{\mathfrak{C}}

\newcommand{\id}{\rm{id}}

\newcommand{\mat}{\mathbf{Mat}}

\newcommand{\ev}{\mathbf{Vect}}
\newcommand{\Ev}{\mathbf{VECT}}

\newcommand{\dev}{{\sf 2GVect}}
\newcommand{\Dev}{{\sf 2GVECT}}
\newcommand{\adcat}{{\sf AdCat}}
\newcommand{\cat}{{\sf Cat}}

\newcommand{\qed}{\hspace*{\fill}$\Box$  \ifmmode \else
    \par\addvspace\topsep\fi}

\newenvironment {proof}{\par\addvspace\topsep\noindent{\it Proof.}
    \ignorespaces }{\qed}

\begin{document}

\title{Generalized 2-vector spaces and general linear 2-groups}  
\author{Josep Elgueta \\ Dept. Matem\`atica Aplicada II \\ Universitat
  Polit\`ecnica de Catalunya \\ email: Josep.Elgueta@upc.edu}
\date{}

\maketitle

\begin{abstract}
In this paper a notion of {\sl generalized 2-vector space} is
introduced which includes Kapranov and Voevodsky 2-vector spaces. Various
kinds of generalized 2-vector spaces are considered
and examples are given. The existence of non free generalized 2-vector spaces
and of generalized 2-vector spaces which are non Karoubian (hence, non
abelian) categories is 
discussed, and it is shown how any generalized 2-vector space can be
identified with a full subcategory of an (abelian) functor category  
with values in the category ${\bf VECT}_K$ of (possibly infinite dimensional)
vector spaces. The 
corresponding general linear 2-groups $\G\Lb(\ev_K[\Cc])$ are
considered. Specifically,  
it is shown that $\G\Lb(\ev_K[\Cc])$ always contains as a (non full)
sub-2-group 
the 2-group ${\sf Equiv}_{{\sf Cat}}(\Cc)$ (hence, for finite categories $\Cc$,
they contain {\sl Weyl sub-2-groups}
analogous to the usual Weyl subgroups of the general linear groups), and
$\G\Lb(\ev_K[\Cc])$ is explicitly computed 
(up to equivalence) in a special case of generalized 2-vector spaces which
include those of Kapranov and Voevodsky. Finally, other important drawbacks of
the notion of generalized 2-vector space, 
besides the fact that it is in general a non Karoubian 
category, are also mentioned at the end of the paper. 
\end{abstract}

\section{Introduction}

\paragraph{} For the
development of 1-{\it dimensional} (i.e., categorical) mathematics, where sets
are sistematically 
replaced by categories, it
would be desirable to have an analog of the usual
linear algebra which has proved so useful in the (0-dimensional)
mathematics of sets. The first logical step in the search of such an analog is
to find a good notion of {\it 
categorical vector space}, more often called a {\sf 2-vector space}. 

The notion of (finite dimensional) 2-vector space over a field $K$
was introduced for the first time by Kapranov and  
Voevodsky \cite{KV94}, motivated by Segal's definition of a conformal field
theory \cite{gS88}. The main point in their definition is to take the category
$\ev_K$ of finite 
dimensional vector spaces over $K$ as analog of the field $K$ and to define
a $\ev_K$-module category as a symmetric monoidal category $\V$ (analog of
the abelian group underlying a vector space) equipped
with an action of $\ev_K$ on it (analog of the multiplication by scalars)
satisfying all the usual axioms of a 
$K$-module up to suitable coherent natural isomorphisms
(cf. \cite{KV94} for 
more details). Then, a 2-vector space over $K$ is, according to these authors,
a ``free $\ev_K$-module category of finite rank'', i.e., a $\ev_K$-module
category equivalent in the appropriate sense to $\ev_K^n$ for some $n\geq 0$
(in particular, the 
underlying symmetric monoidal category is a $K$-linear additive category
equipped with the symmetric monoidal structure induced by the direct
sums). When unpacked, however, 
this definition is quite disappointing due to the
long list of required data and coherence axioms.

A simpler and essentially equivalent definition was given
by Neuchl \cite{mN97}, who defined a  
(finite dimensional) 2-vector space over $K$ as a $K$-linear additive
category $\V$  
which admits a (finite) ``basis of absolutely simple objects'', i.e., a
(finite) subset $\BB=\{X_i\}_{i\in I}$ of absolutely simple objects of $\V$
(by which I mean simple objects whose vector spaces of
endomorphisms are 1-dimensional) such that, for any object $X$ of
$\V$, there exists {\it unique} natural numbers $n_i\geq 0$, $i\in I$, all but
a finite number of them zero, so that $X\cong\oplus_{i\in I} X_i^{n_i}$,
where $X_i^{n_i}$ denotes any biproduct of $n_i$ copies of $X_i$, $i\in
I$. If such a basis $\BB$ exists, it is shown that $\V$ is indeed $K$-linear
equivalent to a category $\ev_K^n$ for some $n\geq 0$.

Another definition of 2-vector space over $K$ was introduced
almost a decade later by Baez and Crans in an attempt to define {\sf Lie
  2-algebras} \cite{BC03}. These authors defined a
2-vector space over $K$ as a category in $\ev_K$, and they proved that
the appropriately defined 2-category of such 2-vector spaces is
biequivalent to the familiar 2-category of length one complexes of vector
spaces over $K$. 

The purpose of this paper is to introduce a generalized notion of 2-vector
space which is in the same spirit as that of Kapranov and Voevodsky
and which includes Kapranov and Voevodsky 2-vector spaces. Thus, instead of
categorifying the usual notion of 
$K$-module, with $K$ replaced by $\ev_K$, we pay attention to the fact that
any vector space is, up to isomorphism, the 
set $K[X]$ of all finite formal linear combinations of elements of some set
$X$ with coefficients in 
$K$, equipped with the obvious sum and multiplication by scalars, and we
categorify such a constructive definition. The starting point now is
going to be not a set $X$ but a category $\Cc$. Then, a {\sl generalized
2-vector space over} $K$ can be defined as a $K$-linear additive category $\V$
which is $K$-linear equivalent 
to the {\it free $K$-linear additive category} generated by $\Cc$, for some
category $\Cc$. By analogy with $K[X]$, such a freely
generated $K$-linear additive category is denoted by $\ev_K[\Cc]$. As shown
below (Proposition~\ref{2-ev_KV}), the $K$-linear additive categories
$\ev_K^n$ ($n\geq 0$) underlying Kapranov and Voevodsky 2-vector spaces 
are recovered (up to $K$-linear equivalence) as the categories $\ev_K[\Cc]$
for $\Cc$ a finite discrete 
category. But not all $K$-linear additive categories of the type $\ev_K[\Cc]$,
for $\Cc$ an 
arbitrary category, are of this type. Thus, it is shown with examples
that, in some cases, there also
exists a basis whose objects, however, are non absolutely
simple. Moreover, for an arbitrary category $\Cc$, it is likely
that there exists no basis in $\ev_K[\Cc]$, 
either of absolutely simple objects or not. Arguments in favour of
this possibility are discussed in the sequel.

Together with the vector space $K[X]$, there is another vector space
that can be built from an arbitrary set $X$. Namely, the vector space $K^X$ of
all functions on $X$ with values in $K$. For finite sets, both
vector spaces are isomorphic (actually, both are functorial and define
functors $K[-],K^{(-)}:{\bf
  FinSets}\To\ev_K$ which are naturally isomorphic). It is then
natural to consider also the analog in the category setting of this
vector spaces of functions, namely, the functor categories with values
in $\ev_K$, and to compare both constructions. In
contrast to what happens with vector spaces, however, they are no
longer equivalent, even if we restrict to finite categories. More
precisely, it will be shown that, for a finite category $\Cc$, $\ev_K[\Cc]$
is equivalent to just a certain (full) subcategory of the functor
category $\Ev_K^{\Cc^{op}}$.  

My motivation for introducing this notion of generalized 2-vector
space was the desire of defining an
analog for 2-groups of the (Frobenius and/or Hopf) group algebras
$K[G]$. Thus, the free vector space construction
$K[X]$ is not just functorial. It actually defines a {\it monoidal} functor
$K[-]:{\bf Sets}\To{\bf VECT}_K$, which moreover is a left adjoint of the
forgetful functor $U:{\bf VECT}_K\To{\bf Sets}$.  The fact that $K[-]$ is
monoidal implies 
that it indeed induces a functor $K[-]:{\bf Monoids}\To{\bf Alg}_K$ between
the category of monoids and that of associative $K$-algebras with unit and
hence, also from the category {\bf Grps} of groups to ${\bf Alg}_K$. In a 
completely analogous way, if $\Dev_K$ denotes the
2-category of the above generalized 2-vector spaces, a monoidal structure on
$\Dev_K$ is explicitly described in \cite{jE6} and it is shown that the
construction $\ev_K[\Cc]$ extends to a
monoidal 2-functor $\ev_K[-]:\cat\To\Dev_K$ which is a left 2-adjoint
of the forgetful 2-functor ${\bf U}:\Dev_K\To\cat$. As in the previous
setting, the fact that 
$\ev_K[-]$ is monoidal implies that for any 2-group $\G$ (more generally,
for any monoidal category), the 2-vector space $\ev_K[\G]$ 
spanned by $\G$ inherits a 2-algebra structure. Thus,
the objects $\ev_K[\G]$ can indeed be considered as analogs of usual group
algebras. 

The previous parallelism between the functor $K[-]$ and
the 2-functor $\ev_K[-]$ and the view of $\ev_K[\G]$ as an analog of usual
group algebras seems enough to make worth exploring this notion of generalized
2-vector space. It also seems worth investigating the
representation theory of 2-groups on these generalized 2-vector spaces, and
to compare the resulting theory with that considered in
\cite{jE4}, where representations on Kapranov and Voevodsky 2-vector spaces
are discussed. As a first step in this direction, the general
linear 2-group $\G\Lb(\ev_K[\Cc])$ (i.e., the 2-group of
$K$-linear selfequivalences of $\ev_K[\Cc]$) is completely computed (up to
equivalence) in the special case 
where $\Cc$ 
is a finite (homogeneous) groupoid $\Gg$. In particular, it is shown that
for any such groupoid the general linear 2-group is always split (for
the definition of split 2-group, see \S 2.8),
generalizing the situation encountered for Kapranov and Voevodsky general
linear 2-groups, which correspond to the case $\Gg$ is a finite discrete
category. The relation between these general linear 2-groups 
$\G\Lb(\ev_K[\Cc])$ and the 2-groups ${\sf Equiv}_{{\sf Cat}}(\Cc)$ of
self-equivalences of $\Cc$ and natural isomorphisms between them is also
discussed, leading naturally to the notion of {\sf Weyl sub-2-group} of
$\G\Lb(\ev_K[\Cc])$ for a finite category $\Cc$, analogous to the Weyl
subgroups of the general linear groups.
 
The notion of generalized 2-vector space, however, has some drawbacks with 
respect to Kapranov and Voevodsky 2-vector spaces. One of them is the fact that
generalized 2-vector spaces are non Karoubian (hence, non abelian) categories
in general. Another one is that they have no dual object in the 
usual sense, except when they are Kapranov and Voevodsky 2-vector
spaces. Finally, the categories of morphisms in $\Dev_K$ between arbitrary
generalized 2-vector
spaces are not always generalized 2-vector spaces. Althought on the one hand
this makes the new notion a quite pathological one, it seems on the other hand
the appropriate notion in order to define a natural analog of group algebras in
the category setting. Another generalization of Kapranov and Voevodsky
2-vector spaces which exhibits a more pleasant behaviour is given by the
Karoubian completion of our generalized 2-vector spaces $\ev_K[\Cc]$. But this
will be discussed elsewhere.   

The outline of the paper is as follows. In Section 2, some facts
concerning $K$-linear additive categories and 2-groups are reviewed (notably,
the classification of 2-groups in terms of suitable 3-cocycles), and a few
elementary results needed later are shown. In Section 3, the
notion of generalized 2-vector space is defined as an
analog for categories of the vector spaces $K[X]$, and 
various kinds of examples are considered. In 
particular, we introduce the notion of {\sl free generalized 2-vector space}
and discuss the possibility that there exists non free generalized 2-vector
spaces. In Section 4, we 
consider the category analog of the vector space $K^X$ of
functions on a set $X$ with values in $K$, namely the functor
categories $\Ev_K^{\Cc^{op}}$, and make explicit the relation with 
the generalized 2-vector space $\ev_K[\Cc]$ generated by $\Cc$. Finally, in
Section 5, the 
general linear 2-group of the generalized 2-vector space generated by a finite
{\it homogeneous} 
groupoid is computed, recovering Kapranov and Voevodsky general linear
2-groups as particular cases. The paper finishes with a few comments on the
relation between our generalized 2-vector spaces and the notion of
$\ev_K$-module category 
and on the above mentioned drawbacks our generalized 2-vector spaces
have with respect to Kapranov and Voevodsky 2-vector spaces.

\section{Preliminaries}

In this section, and unless otherwise indicated, $K$ denotes an arbitrary
commutative ring with unit. 

\paragraph{\S 2.1. $K$-linear additive categories.}

\label{repas_categories_additives}

Recall that a category $\Ll$ is called $K$-{\sf linear} when its 
sets of morphisms come equipped with $K$-module structures such that all
composition maps are $K$-bilinear. When $K=\Integer$, $\Ll$ is often
called a {\sf preadditive category} or an {\sf
Ab-category}.

For any pair of objects $X,Y$ in a $K$-linear category  $\Ll$, a {\sf
biproduct} (or {\sf direct sum}) of $X$ and $Y$ is an object, usually
denoted $X\oplus Y$, together with morphisms $\iota_X:X\To X\oplus Y$,
$\iota_Y:Y\To X\oplus Y$ (called {\sf injections}) and
$\pi_X:X\oplus Y\To 
X$, $\pi_Y:X\oplus Y\To Y$ (called {\sf projections}) such
that   
\begin{equation} \label{equacions_definicio_biproducte}
\pi_X\iota_X={\id}_X,\quad \pi_Y\iota_Y={\id}_Y,\quad
\iota_X\pi_X+\iota_Y\pi_Y={\id}_{X\oplus Y}
\end{equation}
(althought the definition is usually given for preadditive categories, it
actually makes
sense for an arbitrary $K$, the multiplication by scalars playing no role in
the definition). Any diagram in $\Ll$ of the form
$$
\xymatrix{
X\ar@<0.6ex>[rr]^{\iota_X} & & X\oplus
Y\ar@<0.6ex>[ll]^{\pi_X}\ar@<0.6ex>[rr]^{\pi_Y} & &
  Y\ar@<0.6ex>[ll]^{\iota_Y}
}
$$
whose morphisms satisfy (\ref{equacions_definicio_biproducte}) is called a
{\sf biproduct diagram}. The definition extends in the obvious way to
any finite set of objects $X_1,\ldots,X_n$ with $n>2$. 

A {\sf $K$-linear additive category} is a $K$-linear category
$\Ll$ which has a zero object and all binary biproducts (hence, all finite
biproducts). A $\Integer$-linear additive category is usually
called an {\sf additive category}. 

For any finite biproduct $(X_1\oplus\cdots\oplus
X_n,\iota_{X_1},\ldots,\iota_{X_n},
\pi_{X_1},\ldots,\pi_{X_n})$ of $X_1,\ldots,X_n$, 
the tuples $(X_1\oplus\cdots\oplus
X_n,\iota_{X_1},\ldots,\iota_{X_n})$ and
$(X_1\oplus\cdots\oplus
X_n,\pi_{X_1},\ldots,\pi_{X_n})$ 
turn out to be a coproduct
and a product of $X_1,\ldots,X_n$, respectively. By the universal
properties of products and coproducts, this means that the biproduct
of $X_1,\ldots,X_n$ is unique up 
to a unique isomorphism commuting with the injections (or with the
projections). Furthermore, they also make possible to describe a morphism
$f:X_1\oplus\cdots\oplus X_n\To Y_1\oplus\cdots\oplus Y_m$ between biproduct 
objects in terms of an $m\times n$ matrix with entries the composite morphisms
$f_{ij}=\pi_{Y_i}\ f\ \iota_{X_j}:X_j\To Y_i$, $i=1,\ldots,m$,
$j=1,\ldots,n$. Composition is then given by
the formal matrix product and the composition law in $\Ll$. This notation,
however, does not make explicit the injections and projections and must be used
with care.

\paragraph{\S 2.2. The 2-category of $K$-linear additive categories.}
Given $K$-linear categories $\Ll$ and $\Ll'$, a
functor $F:\Ll\To\Ll$ is called $K$-{\sf linear} when it acts $K$-linearly on
the $K$-modules of morphisms. If $K=\Integer$, $F$ is called an {\sf
  additive functor}. It is shown that $K$-linear functors $F:\Ll\To\Ll'$ map
biproduct diagrams to biproduct diagrams and zero objects to
zero objects.

\begin{defn}
Let $\adcat_K$ be the 2-category whose objects and 1- and
2-morphisms are the $K$-linear additive categories, the $K$-linear functors
and all natural transformations, respectively. Composition laws
and identities are the usual ones.
\end{defn}
Observe that $\adcat_K$ is a $K$-{\it linear 2-category}, i.e., all
hom-categories ${\bf Hom}_{\adcat_K}(\Aa,\Aa')$, for $\Aa$ and 
$\Aa'$ any objects in $\adcat_K$, are $K$-linear and the composition
functors are $K$-bilinear. 
  
Among the objects in $\adcat_K$, we have the
category ${\bf Mod}_K$ of all $K$-modules and $K$-linear maps, and the
full subcategory ${\bf Mod}^{f}_K$ of finitely generated
$K$-modules. If $K$ is a field, these categories are denoted
$\Ev_K$ and $\ev_K$, respectively. Objects in $\adcat_K$, for $K$ a
field, further include the categories ${\bf  Rep}_{\ev_K}(G)$ of finite
dimensional linear representations of an arbitrary group $G$.

Observe that, if $\Aa$ and $\Aa'$ are $K$-linear additive categories, the 
product $\Aa\times\Aa'$ inherits a $K$-linear additive structure where
biproducts are given termwise. In particular, the product categories
$\ev_K^r$, for any $n\geq 2$, are also objects in
$\adcat_K$ for $K$ a field. Such objects play an special role in what
follows. 

\paragraph{\S 2.3. Krull-Schmidt $K$-linear additive categories.}
There is a distinguished family of objects in $\adcat_K$ characterized by the
property of having a ``basis''. The formal definition is as follows: 

\begin{defn} \label{definicio_base}
Let $\Aa$ be any object in $\adcat_K$ and $\Sd=\{X_i\}_{i\in
  I}$ any set of objects of $\Aa$. The $K$-{\sf linear
additive subcategory of $\Aa$ generated} (or {\sf spanned}) by $\Sd$ is the
full repletive subcategory of $\Aa$, denoted by $\langle\Sd\rangle$, which
contains a zero 
object {\bf 0} and all biproducts $X_1\oplus\cdots\oplus X_r$ for all objects
$X_1,\ldots,X_r$ in $\Sd$ and all $r\geq 1$ (in particular, if $\Sd=\emptyset$,
$\langle\Sd\rangle$ is a terminal category). When
$\langle\Sd\rangle=\Aa$,  $\Sd$ is said to be an {\sf
additive generating system} or to {\sf additively span} $\Aa$.

A set of objects $\Sd=\{X_i\}_{i\in I}$ of $\Aa$ is called {\sf additively
  free} if, whenever we have an isomorphism $X_{i_1}\oplus\cdots\oplus
  X_{i_r}\cong X_{i'_1}\oplus\cdots\oplus X_{i'_{r'}}$ with
  $X_{i_p},X_{i'_{p'}}\in\Sd$ for all $p=1,\ldots,r$ and $p'=1,\ldots,r'$
  ($r,r'\geq 1$), it
  is $r=r'$ and $X_{i'_{\sigma(p)}}=X_{i_p}$ for some permutation $\sigma\in
  \Sigma_r$ (in particular, the objects in $\Sd$ are non zero and pairwise
  nonisomorphic). 

A (finite) set of objects $\BB=\{X_i\}_{i\in I}$ of $\Aa$ is called a
{\sf (finite) basis} of $\Aa$ if it is additively free and additively spans
$\Aa$ (equivalently, if for any
  nonzero object $X$ there exists unique natural 
  numbers $n_i\geq 0$, $i\in I$, all but a finite number of them zero,
  such that $X\cong \bigoplus_{i\in I}X_i^{n_i}$). When
all objects $X_i$ are simple (resp. simple and with 1-dimensional vector
spaces of endomorphisms), $\BB$ is said to be a {\sf basis of simple objects}
(resp. a {\sf basis of absolutely simple objects}) \footnote{If $\Aa$ is an
abelian category and $K$ is an algebraically closed field, a simple
object is automatically absolutely simple by Schur's lemma. However, this is
not true in general, as it is shown below
(cf. Propositon~\ref{2ev_generat_monoide}).}. 
\end{defn}

The existence of a basis in a $K$-linear additive category is related to a
Krull-Schmidt type theorem. In general, such
theorems have to do with the existence and uniqueness (up to isomorphism and
permutations) of a decomposition as a ``product''
of certain ``indecomposable'' objects of some of the objects in 
certain categories 
(mostly additive categories, but not necessarily). The concrete notions of
product and indecomposable object depend on the particular version of the
theorem. Thus, there is a Krull-Schmidt theorem for the category of
groups in which the product is taken to be
the usual direct product of groups, and where the indecomposable objects are
are the groups $G\neq 1$ such that
$G\cong H\times K$ implies that $H\cong 1$ or $K\cong 1$. The theorem then 
states that any group $G$ satisfying 
either the ascending or descending chain condition on normal subgroups is
isomorphic to a direct product of a finite number of indecomposable groups
and that, when it satisfies both conditions, this decomposition is unique up to
isomorphism and permutation of the factors (see for ex. \cite{tH74}).

For $K$-linear additive categories, one usually takes the biproduct as
the appropriate notion of product, and the
indecomposable objects are the objects $X\ncong{\bf 0}$ such that
$X\cong X'\oplus X''$ implies $X'\cong{\bf 0}$ or 
$X''\cong{\bf 0}$. Clearly, if a basis
$\BB=\{X_i\}_{i\in I}$ indeed exists in such a category, the objects $X_i$
are necessarily indecomposable in this sense (otherwise, it will be $X_i\cong
X\oplus X'$ for some $X,X'\ncong{\bf 0}$ and, hence, $X_i\cong
X_{i_1}\oplus\cdots X_{i_k}\oplus X_{i'_1}\oplus\cdots X_{i'_{k'}}$ for some
$k,k'\geq 1$, in contradiction with the additive freeness of $\BB$). This
suggests introducing the 
following terminology: 

\begin{defn}
A {\sf Krull-Schmidt $K$-linear additive category} is a $K$-linear additive
category which has a basis.
\end{defn}

Notice that any basis $\BB$ in a Krull-Schmidt
$K$-linear additive category $\Aa$ necessarily
includes one (and only one) representative from each isomorphism class of
indecomposable objects (otherwise, $\BB$ will not span $\Aa$
additively). Hence, in contrast to what happens with vector 
spaces, the basis in a Krull-Schmidt $K$-linear 
additive category is unique up to isomorphism. More precisely, we have the
following 

\begin{prop} \label{unicitat_base_2ev}
If $\Aa$ is a Krull-Schmidt $K$-linear additive category, there exists a
unique basis up to isomorphism, given by one representative in each
isomorphism class of indecomposable objects. In 
particular, all bases of $\Aa$ have the same cardinal
(called the {\sf rank of} $\Aa$ and denoted ${\sf rk}(\Aa)$).
\end{prop}

\begin{ex} \label{exemple_2-ev} {\rm 
For any field $K$ and $n\geq 1$, $\ev_K^n$ is a Krull-Schmidt $K$-linear
additive category,
a basis being given by $\BB=\{K(i,n),\ i=1,\ldots,n\}$, with
$K(i,n)=(0,\ldots,\stackrel{i)}{K},\ldots,0)$ for all $i=1,\ldots,n$.
}
\end{ex}

\begin{ex} \label{representacions_de_G} {\rm
If $K$ is an algebraically closed field and $G$ a finite group, ${\bf
  Rep}_{\ev_K}(G)$ is also a Krull-Schmidt $K$-linear additive category, a
  basis being given by one representative in each isomorphism class of simple
  objects, usually called the {\it irreducible representations}
  (there are as many such isomorphism classes as conjugacy classes in $G$;
  see, for ex., Fulton-Harris \cite{FH91}). 
}
\end{ex}

Let us finally point out that Krull-Schmidt $K$-linear additive categories can
also be characterized in terms of the commutative (additive) monoid 
with the isomorphism classes of objects as elements and with the sum
induced by the biproduct of corresponding representative
objects. Specifically, if we denote by ${\bf M}(\Aa)$ this monoid, for any
$K$-linear additive category $\Aa$, we have the following obvious result:

\begin{prop}
A $K$-linear additive category $\Aa$ is of the Krull-Schmidt type if and only
if ${\bf M}(\Aa)$ is free.
\end{prop}
For example, for $n\geq 1$ it is ${\bf M}(\ev_K^n)\cong\Natural^n$, while
  ${\bf M}({\bf Rep}_{\ev_K}(G))\cong\Natural^r$ with 
$r$ the number of conjugacy classes in $G$.

\paragraph{\S 2.4. Free $K$-linear categories.} 

\label{categories_K-lineals_additives_lliures}

For any
category $\Cc$, the {\sf free $K$-linear category} (or {\sf free preadditive
    category} when $K=\Integer$) {\sf generated by} $\Cc$ is the 
$K$-linear category $K[\Cc]$ with the same objects as $\Cc$, with vector
spaces of morphisms 
$$
{\rm Hom}_{K[\Cc]}(X,X'):=K[{\rm Hom}_{\Cc}(X,X')]
$$
and with composition law given by the $K$-bilinear extension of the
composition law in $\Cc$ (identities are the obvious 
ones).

There is an obvious inclusion functor ${\sl k}_{\Cc}:\Cc\To K[\Cc]$, and the
pair $(K[\Cc],{\sl k}_{\Cc})$ 
has the following universal property,  which follows from the
universal property of free $K$-modules: for any $K$-linear 
category $\Ll$ and any functor $F:\Cc\To\Ll$, 
there exists a unique $K$-linear functor $\overline{F}:K[\Cc]\To\Ll$, called
the {\sf $K$-linear extension} of $F$, such that
$F=\overline{F}\ {\sl k}_{\Cc}$. Furthermore, any natural transformation
$\tau:F\Rightarrow G:\Cc\To\Ll$ defines a natural transformation between the
$K$-linear extensions 
$\overline{\tau}:\overline{F}\Rightarrow\overline{G}$.

Note also that the construction $K[\Cc]$ preserves coproducts, i.e., for an
arbitrary family of categories $\{\Cc_i\}_{i\in I}$ it 
is  
\begin{equation} \label{preservacio_coproductes} 
K\left[\sqcup_{i\in I}\Cc_i\right]\ \simeq_K\ \stackrel{K}{\sqcup}_{i\in
  I} K[\Cc_i] 
\end{equation}
where $\simeq_K$ denotes $K$-linear equivalence and $\stackrel{K}{\sqcup}$
denotes the coproduct of $K$-linear categories, given by the usual disjoint
union of categories except that for pairs of objects
in different categories the corresponding hom-set in the coproduct
is the zero vector space, instead of the empty set.

In general, it is possible that non isomorphic objects in $\Cc$ become
isomorphic in $K[\Cc]$. This suggests introducing the following

\begin{defn}
A category $\Cc$ is called $K${\sf -stable} if isomorphic objects in $K[\Cc]$
are also isomorphic in $\Cc$.
\end{defn}

Examples of categories which are $K$-stable for any $K$ include all
groupoids and all free categories. Another example which will be
needed later (see Lemma~\ref{Equiv_vs_Aut}) is provided by the following 

\begin{prop} \label{K-estabilitat}
Let $\Cc$ be a category such that, for any object $X$ of $\Cc$, an
endomorphism $f:X\To X$ is an isomorphism if and only if it is a
monomorphism. Then, $\Cc$ is $K$-stable for any $K$. 
\end{prop}
\begin{proof}
Suppose $X,Y$ are isomorphic objects in $K[\Cc]$, and let
$\sum_{i}\lambda_i f_i:X\To Y$ be an isomorphism, with inverse
$\sum_j\mu_j g_j:Y\To X$. In particular, it is
$$
\sum_{i,j}\lambda_i\mu_j g_jf_i={\id}_X,\quad \sum_{i,j}\lambda_i\mu_j
f_ig_j={\id}_Y
$$
Since the hom-sets in $\Cc$ constitute linear bases
for the corresponding vector spaces of morphisms in $K[\Cc]$, it follows that
there exists pairs $(i_0,j_0)$ 
and $(i_1,j_1)$ such that $g_{j_0}f_{i_0}={\id}_X$ and
$f_{i_1}g_{j_1}={\id}_Y$. In particular, both $f_{i_0}$ and $g_{j_1}$ are
sections (hence, monomorphisms) and consequently, the composite
$f_{i_0}g_{j_1}:X\To X$ is a monomorphism. By hypothesis, $f_{i_0}g_{j_1}$ is
then an isomorphism, from which we conclude that $f_{i_0}$ is an
epimorphism. But a section which is at the same time an epimorphism is
necessarily an isomorphism. Therefore, $X\cong Y$ already in $\Cc$.     
\end{proof}

Let us finally remark that, when the category $\Cc$ is already $K$-linear, the
$K$-linear 
structure on $K[\Cc]$ has nothing to do with that on $\Cc$. Thus, it is a
priori possible that the biproduct of two objects $X,Y$ exists in $\Cc$ while
it does not exist in $K[\Cc]$, and conversely. Similary, there 
is no zero object in $K[\Cc]$ althougt it may exists one in $\Cc$.

\paragraph{\S 2.5. Free additive categories.}
Suppose we are now given a $K$-linear (a preadditive when $K=\Integer$)
category $\Ll$. The 
{\sf free additive category generated by} $\Ll$ is the category
${\sf Add}(\Ll)$ with objects all finite (possibly empty) ordered
sequences of objects in $\Ll$ and with arrows the matrices of arrows in
$\Ll$. More explicitly, a morphism in ${\sf Add}(\Ll)$ between two nonempty
sequences $(X_1,\ldots,X_n)$ and $(X'_1,\ldots,X'_{n'})$ is an $n'\times n$
matrix $A$ whose $(i',i)$-entry is $A_{i'i}\in{\rm
  Hom}_{\Ll}(X_i,X'_{i'})$ (if one or both sequences are empty, the
corresponding hom-set is a singleton, whose element is 
generically denoted by $0$ and called a zero morphism). Composition is  
given by the formal matrix product and the composition law in $\Ll$ when all
involved objects are nonempty, it is equal to the appropriate zero matrix
when only the middle object is empty and it is the corresponding
zero morphim otherwise.

${\sf Add}(\Ll)$ has the obvious
$K$-linear structure inherited from $\Ll$ and the empty sequence as 
a zero object, and it is an additive category, with biproducts given, for
example, by concatenation of sequences. There may exists, however, 
other zero objects (for instance, if $\Ll$ has a zero object {\bf 0}, any
sequence $({\bf 0},\ldots,{\bf 0})$ is also a zero object for ${\sf
  Add}(\Ll)$), and other biproducts (for instance, if $X,X'$ already have a
biproduct $X\oplus X'$ in $\Ll$, a biproduct of $(X)$ and $(X')$ is also given
by the length one sequence $(X\oplus X')$; see
Proposition~\ref{fidelitat_plena_extensio_K-lineal}). 
Note also that any sequence $(X_1,\ldots,X_k)$ of length 
$n\geq 2$ can be thought of as the object part of a biproduct of
the length one sequences $S_1=(X_1),\ldots,S_k=(X_k)$, and that the matrix
$A=(A_{i'i})$ giving a morphism 
$(X_1,\ldots,X_n)\To(X'_1,\ldots,X'_{n'})$ coincides with the matrix
representation of $A$ with respect to the corresponding projections and
injections. 

Two easy facts concerning free additive categories and which will be useful in
the sequel are the following: 

\begin{prop} \label{fidelitat_plena_extensio_K-lineal}
Let $\Ll,\Ll'$ be $K$-linear categories. Then:
\begin{itemize}
\item[(i)] There is a $K$-linear equivalence ${\sf Add}(\Ll)\times{\sf
  Add}(\Ll)\simeq_K{\sf Add}(\Ll\stackrel{K}{\sqcup}\Ll')$, where
  $\stackrel{K}{\sqcup}$ denotes the coproduct of $K$-linear categories
  (see \S 2.4). 
\item[(ii)] For any objects $X,X_1,\ldots,X_n$ in
$\Ll$, the following statements are equivalent:
\begin{itemize}
\item[(a)] $(X)\cong(X_1,\ldots,X_n)$ in ${\sf Add}(\Ll)$.
 \item[(b)] The biproduct of $X_1,\ldots,X_n$ exists in $\Ll$ and 
  $X\cong X_1\oplus\cdots\oplus X_n$.
\end{itemize}
\end{itemize}
\end{prop}
\begin{proof}
A $K$-linear equivalence $E:{\sf Add}(\Ll)\times{\sf
  Add}(\Ll')\To{\sf Add}(\Ll\sqcup\Ll')$  is the functor which maps the object
$((X_1,\ldots,X_n),(X'_1,\ldots,X'_{n'}))$ 
to $(X_1,\ldots,X_n,X'_1,\ldots,X'_{n'})$ and a morphism $(A,A')$ to the
morphism $A\oplus A'$, the usual direct sum of matrices. Such a functor is
indeed essentially surjective because any object $(Y_1,\ldots,Y_m)$ in ${\sf
  Add}(\Ll\sqcup\Ll')$ is isomorphic to any of its permuted sequences. The
proof of (ii) is an easy check left to the reader.
\end{proof}

A feature worth emphasizing of the free additive categories ${\sf Add}(\Ll)$
is that the associated monoid ${\bf M}({\sf 
  Add}(\Ll))$ is not necessarily equal to the free commutative monoid
generated by the isomorphism classes of objects in $\Ll$. Thus, 
for an arbitrary $K$-linear category 
$\Ll$, it may happen that two objects $(X_1,\ldots,X_n)$ and
$(X'_1,\ldots,X'_{n'})$ in ${\sf Add}(\Ll)$ are isomorphic 
even when $n\neq n'$. Examples of this naturally arise when
$\Ll$ is already an additive category, as shown by the previous
Proposition. Actually, it seems to be false that ${\sf Add}(\Ll)$ is always 
a Krull-Schmidt $K$-linear additive category, in spite that it is 
freely generated as an additive category. 

Finally, let ${\sl a}_{\Ll}:\Ll\To{\sf
Add}(\Ll)$ be the $K$-linear embedding mapping $\Ll$ into the 
full subcategory with objects the length one sequences. Then, the pair $({\sf
Add}(\Ll),{\sl a}_{\Ll})$ 
has the following universal property, which justifies the name
given to ${\sf Add}(\Ll)$:  

\begin{prop} \label{propietat_universal_Add(C)}
For any $K$-linear additive category $\Aa$ and any
$K$-linear functor $F:\Ll\To\Aa$, there exists a $K$-linear functor
$\hat{F}:{\sf Add}(\Ll)\To\Aa$, {\it unique up to isomorphism}, such that
$F=\hat{F}\ {\sl a}_{\Ll}$ (we shall call $\hat{F}$ a $K$-{\sf linear
extension of} $F$). Furthermore, given a second $K$-linear functor 
$F':\Ll\To\Aa$, a natural transformation $\tau:F\Rightarrow F'$ and any
$K$-linear extensions $\hat{F}$ and $\hat{F}'$ of $F$ and $F'$, respectively,
there exists a unique natural transformation
$\hat{\tau}:\hat{F}\Rightarrow\hat{F}'$ extending $\tau$ (i.e., such that
$\tau=\hat{\tau}\circ 1_{{\sl a}_{\Ll}}$), and $\hat{\tau}$ is an isomorphism
if $\tau$ is an isomorphism. 
\end{prop}

\paragraph{\S 2.6. Notion of 2-group and the 2-category {\sf 2Grp} of
  2-groups.} 
There are various definitions of (weak) 2-group, depending on the amount of
structure assumed on it. In this paper, by a {\sf 2-group} (also called a {\sf
  categorical group}) we shall mean a
monoidal category $(\G,\otimes,I,a,l,r)$ satisfying the following additional
conditions: (1) $\G$ is a groupoid, and (2) any object $A$ of $\G$ is
invertible, in the sense that the functors $-\otimes A,A\otimes-:\G\To\G$
are equivalences. When the monoidal category is strict (i.e., $a$, $l$ and
$r$ are identities) and any object $A$ is {\it strictly} invertible, in the
sense that the functors 
$-\otimes A,A\otimes-$ are not only equivalences but isomorphisms, the 2-group
is said to be {\sf strict}. For example, if $\Cgg$ is any bicategory and $X$
any object of $\Cgg$, the category 
${\sf Equiv}_{\Cgg}(X)$ with objects the autoequivalences $f:X\To X$ and
morphisms all 2-isomorphisms between these is a 2-group, the composition
and the tensor product being respectively given by the vertical composition of
2-morphisms and the composition of 1-arrows and horizontal composition of
2-arrows, and with ${\id}_X$ as unit object (actually, any 2-group is of this
type for some bicategory $\Cgg$ and some object $X$ of $\Cgg$). In case $\Cgg$
is a strict 2-category, the full subcategory ${\sf Aut}_{\Cgg}(X)$ of ${\sf
Equiv}_{\Cgg}(X)$ with objects only the strict invertible endomorphisms of $X$
is a strict 2-group. 

2-groups are the objects of a 2-category {\sf 2Grp}, whose 1-arrows are 
the monoidal functors between the underlying monoidal categories and whose
2-arrows are the monoidal natural transformations between these (for the
precise definitions, see for instance \cite{jE4}). It is shown that any
2-group is equivalent (in {\sf 2Grp}) to a strict 2-group.

\paragraph{\S 2.7. Classification of 2-groups up to equivalence.}
It is a fundamental result in the theory of 2-groups, first proved aparently
by Sinh \cite{hxSi75}, that these are completely
classified (up to the corresponding notion of equivalence in {\sf 2Grp}) by
triples $(G,M,[\alpha])$, with $G$ a group, $M$ a $G$-module and $[\alpha]\in
H^3(G,M)$. For a given 2-group $\G$, the corresponding group $G$ and
$G$-module $M$ are usually denoted by $\pi_0(\G)$ and $\pi_1(\G)$ and called
the {\sf homotopy groups} of $\G$. They
are respectively equal to the group of isomorphism classes of objects of $\G$
(with the product defined by $[A][B]=[A\otimes B]$) and the group ${\rm
  Aut}_{\G}(I)$ of automorphisms of the unit object (it is an abelian group
with the product given by the composition of automorphisms). A basic
feature of 2-groups is that they are ``homogeneous'', in the sense that
$\pi_1(\G)\cong{\rm Aut}_{\G}(A)$ for {\it any} object $A$ of $\G$. There are
two particularly important such isomorphisms of groups, denoted
$\delta_A$ and $\gamma_A$ and defined by
\begin{equation} \label{definicio_gamma_delta}
\delta_A(u)=r_A\circ ({\id}_A\otimes u)\circ r_A^{-1},\quad  
\gamma_A(u)=l_A\circ (u\otimes{\id}_A)\circ l_A^{-1}
\end{equation}
for all $u\in\pi_1(\G)$. 
In terms of these isomorphisms, the action of $\pi_0(\G)$ on $\pi_1(\G)$ is
given by 
\begin{equation} \label{definicio_accio}
[A]\cdot u=\gamma_A^{-1}(\delta_A(u))
\end{equation}
for any representative $A\in[A]$. Finally, the cohomology class
$[\alpha]\in H^3(\pi_0(\G),\pi_1(\G))$ classifying $\G$ (also called the
{\it Postnikov invariant} of $\G$; see \cite{lB92}) is basically determined
by the associator $a$ of the underlying monoidal category. More explicitly,
let us choose a representative $A_g$ for each class $g=[A_g]\in\pi_0(\G)$, and
for any other object $A'\in g$, choose an isomorphism $\iota_{A'}:A'\To
A_g$, with $\iota_{A_g}={\id}_{A_g}$. Then, a classifying 3-cocycle $\alpha\in
Z^3(\pi_0(\G),\pi_1(\G))$ is
\begin{equation} \label{3-cocicle}
\alpha(g_1,g_2,g_3)=
\gamma^{-1}_{A_{g_1g_2g_3}}(\widetilde{a}_{g_1,g_2,g_3})\in\pi_1(\G) 
\end{equation}
with $\widetilde{a}_{g_1,g_2,g_3}\in{\rm Aut}_{\G}(A_{g_1g_2g_3})$ defined
by 
\begin{align}
\widetilde{a}_{g_1,g_2,g_3}=\iota_{A_{g_1}\otimes A_{g_2g_3}}\circ
({\id}_{A_{g_1}}\otimes\iota_{A_{g_2}\otimes A_{g_3}})\circ
a_{A_{g_1},A_{g_2},A_{g_3}}\circ(\iota_{A_{g_1}\otimes
  A_{g_2}}^{-1}\otimes{\id}_{A_{g_3}})\circ\iota_{A_{g_1g_2}\otimes
  A_{g_3}}^{-1} \nonumber \\ \mbox{} \label{associador_tilde}
\end{align}
As a consequence of the pentagon axiom on $a$, it is seen that the map
$\alpha:\pi_0(\G)^3\To\pi_1(\G)$ defined in this way is indeed a (normalized)
3-cocycle, whose cohomology class turns out to be independent of the chosen
representatives $A_g$ and isomorphisms $\iota_{A'}$. 

\paragraph{\S 2.8. Split 2-groups.}
A particularly simple type of 2-groups are those for which the Postnikov
invariant is trivial, i.e. $[\alpha]=0$. It is easily seen (cf. \cite{jE4})
that these are exactly 
the 2-groups equivalent (in {\sf 2Grp}) to skeletal strict
2-groups, i.e., to strict 2-groups whose underlying categories are skeletal
(isomorphic objects are equal). These
2-groups are called {\sf split} 
because a strict 2-group is of this kind when a certain exact
sequence of 2-groups splits. Specifically, for any group $G$ and any abelian 
group $A$, let $G[0]$ be the group $G$ thought of as a discrete 2-group, and
let $A[1]$ be the group $A$ thought of as a 2-group with only one object and
$A$ as group of automorphisms of the unique object. Then, for any 2-group $\G$
there is an inclusion of 
2-groups $\pi_1(\G)[1]\hookrightarrow\G$ and a ``projection'' morphism
$p:\G\To\pi_0(\G)[0]$, the last one mapping each object $A$ of $\G$ to the
corresponding isomorphism class $[A]$. Together, they define a sequence of
2-group morphisms 
\begin{equation} \label{descomposicio_Postnikov_2grup}
{\bf 1}\To\pi_1(\G)[1]\To\G\To\pi_0(\G][0]\To{\bf 1}
\end{equation}
which is exact in the sense that $\pi_1(\G)[1]$ is equivalent to
the kernel of $p$ (i.e., the homotopy fiber of $p$ over the unit object
$[I]$ of $\pi_0(\G)[0]$). We have then the following:

\begin{prop} \label{2grup_escindit}
Suppose $\G$ is a strict 2-group (as pointed out before, this implies no loss
of generality). Then, if there exists a strict section for the exact sequence
(\ref{descomposicio_Postnikov_2grup}) (i.e., a strict monoidal
functor $S:\pi_0(\G)[0]\To\G$ such that $pS={id}_{\pi_0(\G)[0]}$), $\G$ is
split. 
\end{prop}
\begin{proof}
The existence of such a strict monoidal functor $S$ amounts to
the existence of a choice of representatives $A_g$ compatible
with the tensor product, i.e., such that
$A_{g_1g_2}=A_{g_1}\otimes A_{g_2}$ and $A_e=I$. In this case, it readily
follows from (\ref{definicio_gamma_delta}), (\ref{3-cocicle}),
(\ref{associador_tilde}) and the fact that $\G$ is strict that $\alpha$ maps
each triple $(g_1,g_2,g_3)$ to the identity of $\pi_1(\G)$ and hence,
$[\alpha]=0$.
\end{proof}

Note that for strict 2-groups $\G$, not only $\pi_0(\G)$ but also the set
$|\G|$ of objects of $\G$ inherits a group structure given by the
tensor product. In these cases, the group $\pi_0(\G)$ is nothing but the
quotient of 
$|\G|$ modulo the normal subgroup of objects isomorphic to
the unit object $I$ of $\G$ (see \cite{jE4}, \S 3.7). The existence
of the above strict section $S:\pi_0(\G)[0]\To\G$ then corresponds to the
existence of a group morphism section $s:\pi_0(\G)\To|\G|$.

\section{The 2-category ${\sf 2GVECT}_K$}

\label{2-categoria_2ev}

From now on, $K$ denotes an arbitrary field. 

\paragraph{\S 3.1. Notion of generalized 2-vector space, universal property and
  stability under categorical products.}
For any nonempty category $\Cc$, let
$$
\ev_K[\Cc]:={\sf Add}(K[\Cc])
$$
Thus, an object in $\ev_K[\Cc]$ is any finite (possibly empty) ordered 
sequence $(X_1,\ldots,X_n)$ ($(X_i)_n$ for short) of objects of $\Cc$, with
$n\geq 0$, and a morphism between two nonempty objects $(X_i)_n$ and  
$(X'_{i'})_{n'}$ is an $n'\times n$ matrix
  $A=(A_{i'i})$ with $(i',i)$-entry $A_{i'i}\in K[{\rm
  Hom}_{\Cc}(X_i,X'_{i'})]$, i.e., of the form
\begin{equation} \label{morfisme_generic_Vect[C]}
A_{i'i}=\sum_{\alpha=1}^{d_{i'i}} \lambda(i',i)_{\alpha}\
f_{i'i}^{\alpha},\quad \lambda(i',i)_{\alpha}\in K
\end{equation}
where $f_{i'i}^{\alpha}:X_i\To X'_{i'}$ is a morphism in $\Cc$ for all
$\alpha$. Composition is given by the composition law in $\Cc$ (extended
$K$-bilinearly) and the formal matrix product. Observe that the empty sequence
is the unique zero object of $\ev_K[\Cc]$, because $K[\Cc]$ has no zero
object. Finally, we shall agree that $\ev_K[\emptyset]$ is
the terminal category {\bf 1}.

$\ev_K[\Cc]$ may be thought of as an analog of the vector spaces $K[X]$
constructed from arbitrary sets $X$ (more properly, it is an analog of
$\Natural[X]$). The fact that any vector space is of this 
kind up to isomorphism suggests the following definition of generalized
2-vector space: 

\begin{defn} \label{definicio_2-ev}
A {\sf generalized 2-vector space over} $K$ is a $K$-linear additive
category $\V$ which is $K$-linear equivalent to $\ev_K[\Cc]$ for some
category $\Cc$. When $\Cc$ is finite, $\ev_K[\Cc]$ is called a {\sf finitely
  generated generalized 2-vector space}. In particular, terminal categories
(isomorphic to $\ev_K[\emptyset]$) are called {\sf zero 2-vector spaces}.   
\end{defn}

Generalized 2-vector spaces have the following universal
property, which is an immediate consequence of the universal properties of the
pairs $(K[\Cc],k_{\Cc})$ and $({\sf Add}(K[\Cc]),a_{K[\Cc]})$ and
the analog of the fact that a $K$-linear map between vector
spaces is uniquely determined by the image of a basis:

\begin{prop} \label{propietat_universal_vect_C}
For any category $\Cc$, let $\beta_{\Cc}:\Cc\To\ev_K[\Cc]$ denote the
canonical embedding given by the composite
$\Cc\stackrel{k_{\Cc}}{\To}K[\Cc]\stackrel{a_{K[\Cc]}}{\To}\ev_K[\Cc]$. Then,
for any $K$-linear additive category $\Aa$ and any functor
$F:\Cc\To\Aa$, there exists a
$K$-linear functor $\breve{F}:\ev_K[\Cc]\To\Aa$ (called a {\sf $K$-linear
extension of} $F$), unique up to
isomorphism, such that $F=\breve{F}\ \beta_{\Cc}$. Furthermore, given a second
functor $F':\Cc\To\Aa$, a natural transformation (resp. isomorphism)
$\tau:F\Rightarrow F'$ and 
any $K$-linear extensions $\breve{F}$ and $\breve{F}'$ of $F$ and $F'$,
respectively, there exists a unique 
natural transformation (resp. isomorphism)
$\breve{\tau}:\breve{F}\Rightarrow\breve{F}'$ such that 
$\tau=\breve{\tau}\circ 1_{\beta_{\Cc}}$. 
\end{prop}

It readily follows from this universal property that, for any
category $\Cc$ and any $K$-linear additive category $\Aa$, there exists a
$K$-linear equivalence 
$$
{\bf Hom}_{\adcat_K}(\ev_K[\Cc],\Aa)\simeq_K {\bf
  Hom}_{\cat}(\Cc,\Aa)
$$
Furthermore, it is easily seen that the construction $\ev_K[\Cc]$ extends to a
2-functor 
$\ev_K[-]:\cat\To\adcat_K$ (for details, see \cite{jE6}). Thus, if $\Cc$ and
$\Cc'$ are equivalent 
categories, the corresponding generalized 2-vector spaces $\ev_K[\Cc]$ and
$\ev_K[\Cc']$ 
are also equivalent objects in $\adcat_K$. Hence, it can be assumed
without loss of generality that all involved categories are skeletal.

Note also that generalized 2-vector spaces are stable under finite
products. More explicitly, we have the following

\begin{prop} \label{lema_estabilitat_productes}
For any categories $\Cc_1,\ldots,\Cc_n$ there is a
$K$-linear equivalence
\begin{equation} \label{productes_a_coproductes}
\ev_K[\Cc_1]\times\cdots\times\ev_K[\Cc_n]\simeq_K
\ev_K[\Cc_1\sqcup\cdots\sqcup\Cc_n]
\end{equation}
In particular, the cartesian product of a finite number of generalized
2-vector spaces is a generalized 2-vector space.
\end{prop}
\begin{proof}
It is a direct consequence of (\ref{preservacio_coproductes}) and
Proposition~\ref{fidelitat_plena_extensio_K-lineal}-(i). 
\end{proof}

\paragraph{\S 3.2. Kapranov and Voevodsky 2-vector spaces.}
The simplest examples of generalized 2-vector spaces are those generated
by finite discrete categories. These turn out to be the $K$-linear additive
categories 
$\ev_K^n$ ($n\geq 0$) underlying usual Kapranov and Voevodsky 2-vector
spaces. More generally, for any (non necessarily finite) set $X$ and any
$K$-linear additive category $\Aa$, let
$\Aa^{\oplus X}$ be the full subcategory of $\prod_{x\in X}\Aa$ with
objects the ordered sequences $(A_x)_{x\in X}$ of objects in $\Aa$ such
that $A_x={\bf 0}$ (the zero object) for all but a finite number of
$x\in X$. Then, we have the following 

\begin{prop} \label{2-ev_KV}
For any set $X$, let $X[0]$ denote the corresponding discrete 
category. Then, $\ev_K[X[0]]\simeq_K\ev_K^{\oplus X}$, where $\simeq_K$ means
$K$-linear equivalence. In particular, for a finite set $X$ of cardinal $n\geq
1$, it is $\ev_K[X[0]]\simeq_K\ev_K^n$.
\end{prop}
\begin{proof}
Let $\mat_K$ be the category with objects the natural numbers and
morphisms between non zero objects $n\To m$ the $m\times n$ matrices
with entries in $K$. $\mat_K$ is $K$-linear 
equivalent to $\ev_K$, so that it is enough to see that
$\mat_K^{\oplus X}\simeq_K\ev_K[X[0]]$. Such an
equivalence is defined as follows:
\begin{itemize}
\item map the object $(k_x)_{x\in X}$ of $\mat_K^{\oplus X}$ to the finite
  sequence $(x,\stackrel{k_x)}{\ldots},x)_{x\in X}$ (in
particular, the zero object $(\ldots,0,\ldots)$ is mapped to the empty
sequence);
\item map the morphism
  $(A_x)_{x\in X}:(k_x)_{x\in X}\To(k'_x)_{x\in X}$, with $A_x$ a
  $k'_x\times k_x$ matrix, to the morphism 
  $A:(x,\stackrel{k_x)}{\ldots},x)_{x\in X}
\To(x,\stackrel{k'_x)}{\ldots},x)_{x\in X}$
  given by
$$
A=\left(\begin{array}{ccc} A_{x_1} & \cdots & 0 \\ \vdots & \ddots & \vdots \\ 
0 & \cdots & A_{x_n}\end{array}\right)
$$
where $x_1,\ldots,x_n$ are the elements $x\in X$ for which $k_x,k'_x\neq 0$
(here, the entries in $A_{x_i}$ have to be thought of as the corresponding
scalar multiples of ${\id}_{x_i}$, for each $i=1,\ldots,n$).
\end{itemize}
The functor so defined is clearly fully faithful and it is essentially
surjective because any nonempty sequence 
$(x_{i_1},\ldots,x_{i_r})$ is isomorphic in $\ev_K[X[0]]$ to any one of its
permuted sequences.
\end{proof}

A fundamental feature of the finitely generated 2-vector spaces $\ev_K^n$
($n\geq 1$) is that they have a finite basis of 
absolutely simple objects (see Definition~\ref{definicio_base}), given by
$\BB=\{K(i,n)\}_{i=1}^n$, with   
$K(i,n)=(0,\ldots,\stackrel{i)}{K},\ldots,0),\quad i=1,\ldots,n$ (see
Example~\ref{exemple_2-ev}). In fact, this property characterizes them up to
$K$-linear equivalence. Indeed, any
$K$-linear additive category $\Aa$ having a finite 
(possibly empty) basis of absolutely simple objects turns out to be $K$-linear
equivalent to $\ev_K^n$ for some $n\geq 0$. This suggests introducing the
following 

\begin{defn}
A generalized 2-vector space $\V$ is called an {\sf absolutely simple free
  2-vector space} if it has a basis of absolutely simple
  objects. If, moreover, the basis (which is unique up to isomorphism by
  Proposition~\ref{unicitat_base_2ev}) is finite, $\V$ is called {\sf of
  finite rank} or a {\sf Kapranov and Voevodsky 2-vector space}. 
\end{defn}
 
\begin{ex} {\rm
More examples of generalized 2-vector spaces of the Kapranov and Voevodsky 
type are given by the above mentioned representation categories of finite
groups; see Example~\ref{representacions_de_G}.
}
\end{ex}

\paragraph{\S 3.3. Generalized free 2-vector spaces.} 
Not every generalized 2-vector space is a Kapranov and Voevodsky 
2-vector space. For example, there exists generalized 2-vector spaces which
have a basis but whose basic objects are non 
absolutely simple. Among such generalized 2-vector spaces, we have those
generated by non 
trivial monoids. Specifically, let $M[1]$, for any monoid $M$, be the category
with only one object $\ast$ and the elements of $M$ as morphisms. Then, we
have the following 

\begin{prop} \label{2ev_generat_monoide}
For any non trivial monoid $M$, $\BB=\{(\ast)\}$ is a basis of non
absolutely simple objects for $\ev_K[M[1]]$. 
\end{prop}
\begin{proof}
$\Sd=\{(\ast)\}$ clearly spans additively the category
$\ev_K[M[1]]$, so that we only need to see that it is additively
free, which in this case means to see that $\ev_K[M[1]]$ is skeletal, i.e., 
that $(\ast,\stackrel{k)}{\ldots},\ast)\cong(\ast,\stackrel{k')}{\ldots},\ast)$
implies $k=k'$. To see this, observe that a morphism
$A=(m_{i'i}):(\ast,\stackrel{k)}{\ldots},\ast)
\To(\ast,\stackrel{k')}{\ldots},\ast)$ can be
thought of as a $K[M]$-linear map between the free $K[M]$-modules $K[M]^k$ and
$K[M]^{k'}$, with the composition in $\ev_K[M[1]]$ corresponding to the
composition of linear maps. Hence, the condition 
$(\ast,\stackrel{k)}{\ldots},\ast)\cong(\ast,\stackrel{k')}{\ldots},\ast)$ in
$\ev_K[M[1]]$ is equivalent to the
condition $K[M]^k\cong K[M]^{k'}$ as free $K[M]$-modules. The result follows
then from the fact that the ring $K[M]$ has the invariant dimension property
\footnote{Recall that a ring $R$ is said to have the invariant dimension
  property if for any free $R$-module $F$, any two bases of $F$ have the same
cardinality. \label{nota_peu_pagina}}, and this in turn follows from the
fact that $K[M]$ has a homomorphic image,
namely $K$, which is a division ring (see Hungerford \cite{tH74}, Ch. IV, \S
2). Furthermore, $(\ast)$ is clearly a non absolutely simple object
because it has $K[M]$ as vector space of endomorphisms, which is of
dimension $>1$ for any non trivial 
monoid $M$ (if $M$ is trivial, we recover Kapranov and
Voevodsky 2-vector space $\ev_K$). In fact, $(\ast)$ is neither a
simple object, because any element $a\in K[M]$ which is left cancellable but
not a unit defines a monomorphism $a:(\ast)\To(\ast)$ which is not an
isomorphism. 
\end{proof}

This suggests introducing the following more general notion of a free 2-vector
space:  

\begin{defn}
A generalized 2-vector space $\V$ is called {\sf free} when the underlying
$K$-linear additive category is of the Krull-Schmidt type (i.e., it has a
basis, possibly of non absolutely simple 
objects). If it has a finite basis, it is called {\sf of finite
rank} (equal to the cardinal of any basis). Otherwise, it is called {\sf of
infinite rank}.  
\end{defn}

Notice that, defined in this way, a generalized free 2-vector space may
simultaneously be 
non finitely generated and of finite rank. For instance, if 
$M$ is a non finite monoid, it follows from the previous result that
$\ev_K[M[1]]$ is a non finitely generated free 2-vector space of rank
one.

Generalized free 2-vector spaces, as well as generalized free 2-vector spaces
of finite rank, constitute
a subclass of the class of all generalized 2-vector spaces which remains
stable under finite products, i.e.

\begin{prop} \label{lema} 
If $\V,\V'$ are both (finite rank) free
2-vector spaces, with bases $\BB=\{X_i\}_{i\in I}$ and
$\BB'=\{X'_{i'}\}_{i'\in I'}$, respectively, then $\V\times\V'$ is also a
(finite rank) free 2-vector space, with basis
$\BB\cup\BB'\equiv\{(X_i,0')\}_{i\in 
  I}\cup\{(0,X'_{i'})\}_{i'\in I'}$ ($0$ and $0'$ stand for zero objects in
$\V$ and $\V'$, respectively). 
\end{prop}
\begin{proof}
It is easy check left to the reader. 
\end{proof}

The above examples $\ev_K[M[1]]$ of generalized finite rank free 2-vector
spaces which are not of the Kapranov and Voevodsky type are special cases of
the following

\begin{prop} \label{exemple_general_2ev_feblement_lliure}
For any category $\Cc$ with finitely many isomorphism classes of objects
and such that ${\rm Hom}_{\Cc}(X,X')=\emptyset$ when 
$X\ncong X'$, $\ev_K[\Cc]$ is a generalized finite rank free 2-vector space, a
basis being given by any family of length one 
sequences $\{(X_1),\ldots,(X_n)\}$ with $\{X_1,\ldots,X_n\}$ a set of
representative objects of $\Cc$.
\end{prop}
\begin{proof}
Let $\{X_1,\ldots,X_n\}$ be any set of representative objects in $\Cc$. Since
${\rm Hom}_{\Cc}(X,X')=\emptyset$ when $X\ncong X'$, we have
$\Cc\ \simeq\ \sqcup_{i=1}^n\ M_i[1]$, where $M_i={\rm End}_{\Cc}(X_i)$
($i=1,\ldots,n$). Then, the statement readily follows from 
Proposition~\ref{lema_estabilitat_productes}, \ref{2ev_generat_monoide} and
\ref{lema}. 
\end{proof}

Finally, observe that, besides the universal property in
Proposition~\ref{propietat_universal_vect_C}, generalized free 2-vector spaces 
satisfy the following additional universal property:

\begin{prop} \label{propietat_universal_2-ev_lliures}
Let $\V$ be a generalized free 2-vector space, with basis $\BB$, and let
$\V_{\BB}$ be the 
full subcategory of $\V$ with $\BB$ as set of objects. Then,
for any $K$-linear additive category $\Aa$ and any $K$-linear functor
$F:\V_{\BB}\To\Aa$, there exists a
$K$-linear functor $\widetilde{F}:\ev_K[\Cc]\To\Aa$ (called also a {\sf
  $K$-linear extension of} $F$), unique up to
isomorphism, such that $F=\widetilde{F}\ \beta_{\Cc}$. Furthermore, given a
second $K$-linear 
functor $F':\V_{\BB}\To\Aa$, a natural transformation 
$\tau:F\Rightarrow F'$ and 
any $K$-linear extensions $\widetilde{F}$ and $\widetilde{F}'$ of $F$ and $F'$,
respectively, $\tau$ extends uniquely to a
natural transformation
$\widetilde{\tau}:\widetilde{F}\Rightarrow\widetilde{F}'$. 
\end{prop}
\begin{proof}
It is left to the reader.
\end{proof}

\paragraph{\S 3.4. On the existence of generalized non-free 2-vector spaces.}
At this point, the question naturally arises whether {\it any} generalized
2-vector space is free \footnote{Conversely, we can ask
whether any Krull-Schmidt $K$-linear additive category is a generalized
2-vector space.}.   

$\ev_K[\Cc]$ always has an additive
generating system, whatever the category $\Cc$ is. For instance, the set 
of all length one sequences. Moreover, since the empty sequence is
the unique 
zero object in $\ev_K[\Cc]$, all indecomposable sequences $S$
are of length one. 
Hence, determining if it indeed exists a basis in $\ev_K[\Cc]$ and finding it
explicitly involves in the following two steps: (1) to determine which
length one sequences are indecomposable, and (2) to see that the 
indecomposable length one sequences are indeed additively free. 

As regards the first step, note that in all examples of generalized free
2-vector spaces considered until now (as well as in those considered in \S
3.5), all length one sequences turn out to be indecomposable, so that
a basis is given by any representative set of objects of $\Cc$ (see
Proposition~\ref{unicitat_base_2ev}).

The following result gives sufficient conditions on
$\Cc$, different from those in
Proposition~\ref{exemple_general_2ev_feblement_lliure} and
Proposition~\ref{2ev_generat_2-grup} below, which ensure that all length one
sequences are also indecomposable: 

\begin{prop} \label{cas_C_finita}
Let $\Cc$ be a category all whose hom-sets are finite (in particular, $\Cc$
may be a finite category, but not necessarily) and such that all
monomorphisms $f:X\To X$ are isomorphisms for any object $X$ in $\Cc$. Then,
all length one sequences $(X)$ are indecomposable
objects of $\ev_K[\Cc]$. Consequently, if $\ev_K[\Cc]$ has a 
basis for such a category $\Cc$, it is necessarily given by any family of
length one sequences 
$\BB=\{(X_i)\}_{i\in I}$ with $\{X_i\}_{i\in I}$ a set of
representative objects of $\Cc$.
\end{prop}
\begin{proof}
Let us first see that, for any object $X$ of $\Cc$, $X$ is not a
biproduct in $K[\Cc]$ of objects all of them
nonisomorphic to $X$ (either in $K[\Cc]$ or in $\Cc$, because the isomorphism
classes of objects are the same in both categories by
Proposition~\ref{K-estabilitat}). Indeed, suppose $X$ is a biproduct (in
$K[\Cc]$) of 
$X_1,\ldots,X_n$, and let $\iota_k:X_k\To X$ and $\pi_k:X\To X_k$ be the
corresponding injections and projections, with
$\iota_k=\sum_i\lambda_{ki}f_{ki}$ and $\pi_k=\sum_j\mu_{kj}g_{kj}$, where 
$f_{ki}:X_k\To X$ and $g_{kj}:X\To X_k$ are morphisms in $\Cc$ and
$\lambda_{ki},\mu_{kj}\in K$. They are such that
\begin{align} 
\pi_k\iota_k={\id}_{X_k},&\quad k=1,\ldots,n\label{eq1} \\
\sum_{k=1}^n\iota_k\pi_k&={\id}_X \label{eq2}
\end{align}
Then, it follows from (\ref{eq1}) that
$\sum_{i,j}\lambda_{ki}\mu_{kj}(g_{kj}f_{ki})={\id}_{X_k}$ for all
$k=1,\ldots,n$. But ${\rm
  End}_{K[\Cc]}(X_k)$ is the vector space with basis ${\rm 
  End}_{\Cc}(X_k)$ and hence, for each $k=1,\ldots,n$, there exists at least
one pair $(i_k,j_k)$ such 
that $g_{kj_k}f_{ki_k}={\id}_{X_k}$. In particular,
$f_{ki_k}$ is a section. Similarly, it follows from
(\ref{eq2}) that $\sum_{i,j,k}\lambda_{ki}\mu_{kj}(f_{ki}g_{kj})={\id}_{X}$
and hence, there exists, for at least one value of $k$, at least one
pair $(i'_{k},j'_{k})$ such that
$f_{ki'_{k}}g_{kj'_{k}}={\id}_X$. In particular, 
$g_{kj'_k}$ is also a section (hence, a monomorphism). The argument is now the
same we have made to prove Proposition~\ref{K-estabilitat}. Namely,
the composite 
$f_{ki_k}g_{kj'_k}:X\To X$ is a monomorphism and consequently, an
isomorphism by hypothesis. This implies that
$f_{ki_k}:X_k\To X$, which is a section, is
also an epimorphism and hence, an isomorphism. Therefore, at
least one factor $X_k$ is isomorphic to $X$ in $\Cc$. 

Suppose now that $(X)$ is decomposable in $\ev_K[\Cc]$, i.e., $(X)=S\oplus S'$
for some sequences $S,S'$ in $\ev_K[\Cc]$ both of length $\geq 1$. This means
that there exists objects $X_0,\ldots,X_k$ in $\Cc$, with $k\geq 1$, such that
$(X)\cong(X_0,X_1,\ldots,X_k)$ in 
$\ev_K[\Cc]$. According to
Proposition~\ref{fidelitat_plena_extensio_K-lineal}-(ii), however, this holds
if and only if the biproduct of $X_0,X_1,\ldots,X_k$ exists in
$K[\Cc]$ and $X_0\oplus X_1\oplus\cdots\oplus X_k\cong X$ (in $K[\Cc]$). It
follows then from the previous observation that $X\cong X_i$ for at least one
$i\in\{0,1,\ldots,k\}$. Let us assume that $X=X_0$. Then, for any other
object $Y$ in $\Cc$, we have a linear isomorphism 
$$
{\rm Hom}_{K[\Cc]}(X,Y)\cong{\rm Hom}_{K[\Cc]}(X\oplus X_1\oplus\cdots\oplus
X_k,Y) 
$$
i.e.,
$$
K[{\rm Hom}_{\Cc}(X,Y)]\cong K[{\rm
  Hom}_{\Cc}(X,Y)]\oplus\left(\bigoplus_{i=1}^k K[{\rm
    Hom}_{\Cc}(X_i,Y)]\right)
$$
Since all involved hom-sets are finite, this is an isomorphism of finite
dimensional vector spaces. Hence
$$
{\sf dim}_K\left(\bigoplus_{i=1}^k K[{\rm
    Hom}_{\Cc}(X_i,Y)]\right)=0
$$
for any object $Y$, which implies $k=0$, in contradiction with the fact
that $k\geq 1$. 
\end{proof}

\begin{ex} {\rm
Take $\Cc=\mat_{\F_q}$, the category of matrices with entries in the finite
field $\F_q$ of $q$ elements. This category satisfies none of the conditions
stated in Propositions~\ref{exemple_general_2ev_feblement_lliure} or
\ref{2ev_generat_2-grup}. However, all its hom-sets are finite, and
an endomorphism $A:n\To n$ is a monomorphism if and only if it is an 
isomorphism. Hence, by
Proposition~\ref{cas_C_finita}, all length one sequences $(n)$, with $n\geq 1$,
are indecomposable objects in $\ev_K[\mat_{\F_q}]$. Note that this is not in
contradiction with the fact that $n=1\oplus\stackrel{n)}{\cdots}\oplus 1$ in
$\mat_{\F_q}$, because $(n)\cong(1,\stackrel{n)}{\ldots},1)$ is
equivalent to this equality in $K[\mat_{\F_q}]$, not in $\mat_{\F_q}$.
}
\end{ex}

In general, however, it is false that all length one sequences are
indecomposable. For instance, this is not the case if $\Cc$ is already
additive (see
Proposition~\ref{fidelitat_plena_extensio_K-lineal}-(ii)). Furthermore, we
have already pointed out
in \S 2.5 that, for an arbitrary $K$-linear category $\Ll$, the monoid
${\bf M}({\sf   Add}(\Ll))$ is not always isomorphic to the free commutative
monoid generated by the isomorphism classes of objects in $\Ll$.

As regards the second step above, notice that it is equivalent to the
essential uniqueness part in a Krull-Schmidt theorem for these kind of
$K$-linear additive categories.

There are several versions of this theorem, concerning various types of
$K$-linear 
additive categories. The classical version, which goes back to Schmidt (1913)
and Krull (1925), refers to the abelian categories of modules over a
commutative ring with 
unit $K$, and it states that any $K$-module of finite length (more
generally, any $K$-module which is a direct 
sum of $K$-modules with local endomorphism rings) decomposes as a direct sum of
finitely many indecomposables and that the 
decomposition is unique up to isomorphism and permutation of the direct
summands (see for ex. \cite{sW03} or \cite{lR88}). The result was later
shown for the categories of sheaves by Atiyah \cite{mA56}. More
generally, he proved that in
any {\it exact} category satisfying a suitable 
finiteness condition (called the ``bichain condition''), {\it each} object
has an essentially unique decomposition as a finite direct sum of
indecomposables. A third version can be found in \cite{hB68} (p. 20), where
the essential uniqueness of the decomposition for objects analogous to
the above is demonstrated for any {\it Karoubian} additive category (i.e., any
additive category where all idempotents split).

Althought our categories $\ev_K[\Cc]$ are $K$-linear additive, in general they
are neither Karoubian nor exact (see \S 3.6), so that none of the above three
versions applies. Moreover, the proofs of the essential uniqueness
in these three versions make essential use of the fact that
the endomorphism rings of the involved indecomposable objects are local, while
our endomorphism 
rings ${\rm End}_{\ev_K[\Cc]}(X)=K[{\rm End}_{\Cc}(X)]$ need not be
local, even for finite categories $\Cc$.

This suggests that our generalized 2-vector
spaces $\ev_K[\Cc]$ will have no basis in general, even if we restrict to
finite categories $\Cc$ (with more than one object, of course), and that two
notions of freeness should be distinguished in the category
setting. A notion of ``external'' freeness, when the category is a free object
in the appropriate 2-category, and a
notion of ``internal'' freeness, when there exists in the category a family of
basic objects from which all 
objects can be generated in an essentially unique way with the help of some
given associative operation. An externally free category may be internally non
free. In the context of additive categories, this corresponds to free
additive categories ${\sf Add}(\Ll)$ whose monoid ${\bf M}({\sf Add}(\Ll))$ is
non free. The above mentioned existence
of non finitely generated free generalized 2-vector spaces of finite
rank also fits naturally with this situation.

\paragraph{\S 3.5. More examples of generalized free 2-vector spaces of
  infinite rank.} 
The simplest examples of generalized free 2-vector spaces of infinite rank are
those generated by non finite discrete categories (see
Proposition~\ref{2-ev_KV}). These examples are a special case of a more
general situation. 

Indeed, according to Proposition~\ref{exemple_general_2ev_feblement_lliure},
for any category $\Cc$ equivalent to a finite disjoint union of monoids
(viewed as one object categories), $\ev_K[\Cc]$ is a generalized finite rank
free 2-vector. It turns out that, when the
finite hypothesis is replaced by the assumption that all involved monoids are
isomorphic to the same monoid $M$, the result remains true except that the
generalized 2-vector space is now of a possibly infinite rank. Explicitly, let
us call a category $\Cc$ {\sf homogeneous} when it is a disjoint union of
copies (possibly infinite in number) of the same monoid $M$, and call the
monoid $M$ the {\sf underlying monoid} of $\Cc$. Then, we have the following

\begin{prop} \label{2ev_generat_2-grup}
For any homogeneous category $\Cc$,
$\ev_K[\Cc]$ is a generalized free 2-vector space, a 
basis being 
given by any family of length one sequences $\BB=\{(X_i)\}_{i\in I}$ with
$\{X_i\}_{i\in I}$ a set of representative objects of $\Cc$.
\end{prop}
\begin{proof}
We only need to see
that $\BB$ is additively free. Let $M$ be the underlying monoid of $\Cc$. Then,
the same argument used to prove Proposition~\ref{2ev_generat_monoide} can now
be used to show that if $(X_{i_1},\ldots,X_{i_k})$ and
$(X_{i'_1},\ldots,X_{i'_{k'}})$ are any isomorphic objects in $\ev_K[\Cc]$,
then $k=k'$ (any morphism $A$
between them is thought of as a morphism of free $K[M]$-modules
$K[M]^k\To K[M]^{k'}$). It remains to see that both sequences are 
the same up to isomorphism and permutation of the objects. For any
$j=1,\ldots,k$, let $r_j$ and $r'_j$ be the number of 
copies of $X_{i_j}$ present (up to isomorphism) in the first and second
sequences, respectively. By 
definition, it is $r_j\geq 1$, and since $A$ is invertible, it is also
$r'_j\geq 1$ (otherwise, the corresponding matrix column in $A$ will be
entirely made of zeros and hence, non invertible). By symmetry, we conclude
that both sequences, if isomorphic, necessarily 
contain the same objects (up to isomorphism), but probably with different
multiplicities. Let $X_1,\ldots,X_s$ be the pairwise non isomorphic objects
present in both sequences, and let $p_1,\ldots,p_s$ and $p'_1,\ldots,p'_s$
the respective multiplicities, so that 
$$
(X_{i_1},\ldots,X_{i_k})\cong
(X_1,\stackrel{p_1)}{\ldots},X_1,\ldots,X_s,\stackrel{p_s)}{\ldots},X_s)
$$
and
$$
(X_{i'_1},\ldots,X_{i'_{k'}})\cong
(X_1,\stackrel{p'_1)}{\ldots},X_1,\ldots,X_s,\stackrel{p'_s)}{\ldots},X_s)
$$
Then, the isomorphism $A$ is necessarily of the form
$$
A=\left(\begin{array}{ccc} A_1 & \cdots & 0 \\ \vdots & \ddots & \vdots \\ 
0 & \cdots & A_s\end{array}\right)
$$
with $A_l$ a $p'_l\times p_l$ matrix with entries in $K[M]$,
$l=1,\ldots,s$. Furthermore, since $A$ is an isomorphism, there exists a second
matrix 
$$
A'=\left(\begin{array}{ccc} A'_1 & \cdots & 0 \\ \vdots & \ddots & \vdots \\ 
0 & \cdots & A'_s\end{array}\right)
$$
with $A'_l$ a $p_l\times p'_l$ matrix with entries also in $K[M]$, for all
$l=1,\ldots,s$, such that 
$$
A_lA'_l={\rm Id}_{p'_l},\quad A'_lA_l={\rm Id}_{p_l} \qquad l=1,\ldots,s
$$
Using again the invariance dimension property of $K[M]$ \footnote{See footnote
  \ref{nota_peu_pagina}.}, we conclude that
$p_l=p'_l$ for all $l=1,\ldots,s$. 
\end{proof}

\begin{cor}
For any homogeneous groupoid $\Gg$ (in particular, for any discrete category
or any 2-group), finite or not, $\ev_K[\Gg]$ is a generalized free 2-vector
space, a basis being  
given by any family of length one sequences $\BB=\{(X_i)\}_{i\in I}$ with
$\{X_i\}_{i\in I}$ a set of representative objects of $\Gg$.
\end{cor}

\paragraph{\S 3.6. Existence of generalized 2-vector spaces which are non
  abelian (even non Karoubian) categories.}  
Kapranov and Voevodsky 2-vector spaces are abelian categories (in particular,
Karoubian), but this is no longer true for generalized 2-vector spaces.

To see this, let us consider the case of the generalized 
2-vector spaces $\ev_K[M[1]]$ generated by non trivial monoids (see \S
3.3 for the notation). Given
a non zero endomorphism $a:(\ast)\To(\ast)$, with $a\in K[M]$, we want to
know if $a$ has a kernel or not. Note first the following

\begin{lem} \label{lema_monomorfismes}
Let $K$ be of characteristic $\neq 2$. Then, any monomorphism $S\To(\ast)$
in $\ev_K[M[1]]$, with $S\neq\emptyset$, is necessarily an endomorphism
$b:(\ast)\To(\ast)$ for some $b\in K[M]$.
\end{lem}
\begin{proof}
Let us suppose that there exists a monomorphism 
$A:(\ast,\stackrel{r)}{\ldots},\ast)\To(\ast)$ with $r>1$, given by a row
matrix $A=(a_1\ \cdots\ a_r)$, $a_i\in K[M]$. At least one of the entries
is non zero (otherwise, it will not be monic). Let 
$a_1\neq 0$, so that $a_1\neq -a_1$ ($K$ is of characteristic $\neq 2$). Then,
the morphisms $B,B':(\ast)\To(\ast,\stackrel{r)}{\ldots},\ast)$ given
by
$$
B=\left(\begin{array}{c} a_2 \\ -a_1 \\ 0 \\ \vdots \\ 0\end{array}\right)
\qquad B'=\left(\begin{array}{c} -a_2 \\ a_1 \\ 0 \\ \vdots \\
    0\end{array}\right) 
$$
clearly satisfy $B\neq B'$ and $AB=AB'$, in contradiction with the hypothesis
that $A$ is monic.
\end{proof}

The answer to the above question reads then as follows:

\begin{prop}
Let $M$ be a non trivial monoid and let $a\in K[M]$. Then, as a morphism
$a:(\ast)\To(\ast)$ in $\ev_K[M[1]]$, it has a kernel if and only if $a$ is
not a right zero divisor of $K[M]$, in which case the kernel is the morphism
$\emptyset\To(\ast)$.
\end{prop}
\begin{proof}
If $a$ is not a right zero divisor, the only morphisms
$B:S\To(\ast)$ such that $aB=0$ are the zero morphisms, so 
that $\emptyset\To(\ast)$ is clearly a kernel of $a$.

Suppose now that $a$ is a right zero
divisor. In this case, the zero morphism $\emptyset\To(\ast)$ can no
longer be a kernel of $a$, because there are non zero morphisms $B:S\To(\ast)$
such that $aB=0$. For example, all morphisms $b:(\ast)\To(\ast)$ with $b\in
K[M]$ such that $ab=0$. By the previous Lemma, if a kernel exists, being
monic, it is
necessarily an endomorphism $b:(\ast)\To(\ast)$ for some $b\in K[M]$ such that
$ab=0$. But the universal property of the kernel further 
requires $b$ to be left cancellable, in contradiction with the fact
that it is a left zero divisor.
\end{proof}

For example, if $M=\Integer_2=\{\pm\}$, it is easily checked that  
$(+)-(-)$ is a zero divisor in $K[\Integer_2]$ and hence it has no kernel as
endomorphism of $(\ast)$ in $\ev_K[\Integer_2[1]]$.

It might be thought that the categories $\ev_K[\Cc]$, althought non
abelian in general, they are at least Karoubian. But this is also false, as
Example~\ref{exemple_M=N_continuacio} below shows.

\paragraph{\S 3.7. The 2-category of generalized 2-vector spaces and some full
  sub-2-categories.} 
Let the {\sf 2-category of generalized 2-vector spaces over} $K$, denoted by
$\Dev_K$, be the full sub-2-category of $\adcat_K$ with objects all
generalized 2-vector spaces over $K$. In particular, 1-morphisms in ${\sf
  2GVECT}_K$ are 
$K$-linear functors and 2-morphisms arbitrary natural transformations between
these. As full sub-2-category of $\adcat_K$, observe that $\Dev_K$ is a
$K$-linear 2-category (see \S 2.2).

There are various full sub-2-categories of ${\sf 2GVECT}_K$ that can be
distinghuished, according to the various types of generalized 2-vector spaces
considered before. Thus, let 

\begin{itemize}
\item $\Dev_K^f$ be the full sub-2-category of $\Dev_K$ 
with objects only the generalized free 2-vector spaces;
\item $\Dev_K^{ff}$ be the full sub-2-category of $\Dev_K$ 
with objects only the generalized finite rank free 2-vector spaces;
\item $\dev_K$ be the full sub-2-category of $\Dev_K$ 
with objects only the finitely generated generalized 2-vector spaces;
\item $\dev_K^f$ be the full sub-2-category of $\Dev_K$ 
with objects only the finitely generated generalized free 2-vector spaces
(hence, of a necessarily finite rank), and
\item ${\sf 2Vect}_K$ be the full sub-2-category of $\Dev_K$ 
with objects only the Kapranov and Voevodsky 2-vector spaces.
\end{itemize}
All these sub-2-categories fit into the following diagram of inclusion
2-functors, where the 
label $(\star)$ denotes a strict inclusion while $({\sf id?})$ means that it
could be an identity:
$$
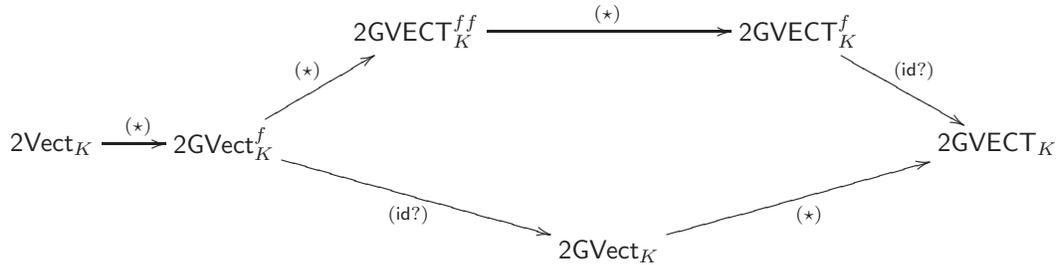
\xymatrix{
& & \Dev_K^{ff}\ar[rr]^{(\star)} & & \Dev_K^f\ar[rd]^{({\sf id?})} & 
\\ {\sf 2Vect}_K\ar[r]^{(\star)} & \dev_K^f\ar[ru]^{(\star)}\ar[rrd]_{({\sf
    id?})} & & & & \Dev_K 
\\ & & & \dev_K\ar[rru]_{(\star)} & & 
} 
$$
Thus, an example of an object in $\Dev_K^f$ which is not in $\Dev_K^{ff}$ is
$\ev_K[X[0]]$ for any infinite set, an example of an object in $\Dev_K^{ff}$
which is not in $\dev_K^f$ is $\ev_K[M[1]]$ for any infinite monoid (or
$\ev_K[\G]$ for any 2-group $\G$ such that $\pi_0(\G)$ is finite and
$\pi_1(\G)$ is infinite), and an example of an object in $\dev_K^f$ which is
not in 
${\sf 2Vect}_K$ is $\ev_K[M[1]]$ for any finite monoid (or $\ev_K[\G]$ for any
finite 2-group). 

\paragraph{\S 3.8. Finite rank free generalized 2-vector spaces up to
  equivalence.} 
An arbitrary free generalized 2-vector space of finite rank encodes a more
involved  
structure than a Kapranov and Voevodsky 2-vector space. Thus,
these 2-vector spaces are completely characterized (up to
equivalence) by their rank (as it occurs for finite dimensional vector
spaces). However, characterizing an arbitrary free generalized 2-vector
space of finite rank (up to equivalence) generally requires a whole set of
structure constants in the field $K$, taking 
account of the non trivial composition law for morphisms between basic
objects. More explicitly, if $\V$ is a free generalized 2-vector space over
$K$ of finite rank $r$ and $\BB=\{X_1,\ldots,X_r\}$ is a basis of
$\V$, we may choose for each pair of basic objects $X_i,X_j\in\BB$ a linear
basis $\BB(X_i,X_j)=\{f(i,j)_{\alpha},\ \alpha\in\Lambda(i,j)\}$ in
the vector space ${\rm Hom}_{\V}(X_i,X_j)$, which we shall assume such
that $f(i,i)_0={\id}_{X_i}$
for all $i\in I$. Then, 
the composition law in $\V$ is completely given by the set of {\sf structure
  constants}  
$\{c(i,j,k)_{\alpha\beta}^{\gamma}\in K,\ i,j,k\in \{1,\ldots,r\},\ 
\alpha\in\Lambda(i,j),\ \beta\in\Lambda(j,k),\ \gamma\in\Lambda(i,k)\}$
defined by
$$
f(j,k)_{\beta}f(i,j)_{\alpha}=
\sum_{\gamma\in\Lambda(i,k)}c(i,j,k)_{\alpha\beta}^{\gamma}f(i,k)_{\gamma} 
$$
These constants satisfy the following associativity and unit equations coming
from the corresponding axioms on the composition law:
\begin{itemize}
\item ({\it associativity}) For all $i,j,k,l\in \{1,\ldots,r\}$,
  $\alpha\in\Lambda(i,j)$, 
  $\beta\in\Lambda(j,k)$, $\gamma\in\Lambda(k,l)$ and
  $\delta\in\Lambda(i,l)$ it is
$$
\sum_{\mu\in\Lambda(i,k)}c(i,j,k)_{\alpha\beta}^{\mu}
c(i,k,l)_{\mu\gamma}^{\delta} 
=\sum_{\nu\in\Lambda(j,l)}c(i,j,l)_{\alpha\nu}^{\delta}
c(j,k,l)_{\beta\gamma}^{\nu}
$$
\item ({\it unit conditions}) For all $i,j\in \{1,\ldots,r\}$ and
  $\alpha,\beta\in\Lambda(i,j)$ it is
$$
  c(i,i,j)_{0\alpha}^{\beta}=\delta_{\alpha\beta},\quad
  c(i,j,j)_{\alpha 0}^{\beta}=\delta_{\alpha\beta}
$$
\end{itemize}
For free generalized 2-vector spaces of rank one, these are nothing but the
equations satisfied by the structure constants of an associative
$K$-algebra with unit.

Althought the above constants depend on the basis $\BB$ of $\V$ and on
the chosen linear basis of morphisms between basic objects, they serve to
completely determine $\V$ in the following sense: 

\begin{prop}
Two finite rank free generalized 2-vector spaces $\V$ and $\V'$ are equivalent
  (as objects in ${\sf 2GVECT}^{ff}_K$) if and only if they have the same
structure constants for suitably chosen bases of objects and linear
  bases of morphisms between basic objects.   
\end{prop}
\begin{proof}
Left to the reader.
\end{proof}

It is worth pointing out, however, that not all sets of constants
satisfying the above associativity and unit conditions are the
structure constants of a finite rank free 2-vector space. For
instance, in the rank one case, it should further exist a
linear basis of endomorphisms of the basic object for which the constants are
given by $c(1,1,1)_{\alpha\beta}^{\gamma}=\delta_{\gamma,m(\alpha,\beta)}$. In 
other words, among all possible associative
algebras with unit, only the algebras of a monoid correspond to free
generalized 2-vector spaces of rank one.

\section{$\ev_K[\Cc]$ versus the functor category
  $\Ev_K^{\Cc^{op}}$.}

Together with $K[X]$, there is one more vector space
which can be built from a set $X$. Namely, the vector space $K^X$ of all
functions on 
$X$ with values in $K$. This construction is also functorial. In fact,
if restricted to finite sets, both functors $K[-],K^{(-)}:{\bf
  FinSets}\To\ev_K$ are naturally isomorphic.

The purpose of this section is to consider the analog for
categories of the 
vector spaces $K^{(X)}$ and to show that the corresponding
construction is no longer equivalent to $\ev_K[\Cc]$, even if we
restrict to finite categories. 

\paragraph{\S 4.1. The functor categories $\Ev_K^{\Cc^{op}}$.} 
To define the analog of $K^{(X)}$, we follow again Kapranov and
Voevodsky insight of replacing $K$ by the category of vector spaces,
except that we shall consider the category ${\bf VECT}_K$ of all vector spaces,
finite dimensional or not. If we further replace the set
$X$ by a category $\Cc$, we are led to the category ${\bf VECT}_K^{\Cc}$ with
objects 
all functors $F:\Cc\To\Ev_K$ and the natural transformations between these as
morphisms. For various reasons, however, it is more convenient to consider the 
category $\Ev_K^{\Cc^{op}}$ of contravariant functors.

Note that $\Ev_K^{\Cc^{op}}$ is a $K$-linear additive
category for any $\Cc$. A zero object is given by
the constant functor $F_0$ mapping each object of $\Cc$ to the zero vector
space, and a biproduct of two functors
$F,G:\Cc^{op}\To\Ev_K$ is given by the composite functor
$\Cc^{op}\stackrel{\Delta}{\To}\Cc^{op}\times\Cc^{op}\stackrel{F\times
  G}{\To}\Ev_K\times\Ev_K\stackrel{\oplus}{\To}\Ev_K$ 
and the obvious injections and projections (here, $\Delta$
denotes the diagonal functor and $\oplus$ the usual direct sum functor on
$\Ev_K$). Such a biproduct is denoted by $F\oplus G$.

Actually, in contrast to the 2-vector spaces $\ev_K[\Cc]$, the categories
$\Ev_K^{\Cc^{op}}$ are always abelian, for any category $\Cc$, because the
target category $\Ev_K$ is already abelian. 

\begin{ex} \label{exemple_M=N} {\rm 
If $\Cc=\Natural[1]$, with $\Natural$ the additive monoid of natural numbers, a
functor $F:\Natural[1]^{op}\To\Ev_K$ is completely given by a vector space $V$
(the image of $\ast$) together with a $K$-linear map $f:V\To V$ (the image of
the morphism $1:\ast\To\ast$), and both can be chosen arbitrarily because
$\Natural$ is free. By identifying $f$ with the action of an
indeterminate $T$ on $V$ and extending this action in the obvious way to the
whole polynomial algebra $K[T]$, objects of
$\Ev_K^{\Natural[1]^{op}}$ are naturally identified with modules over the
polynomial algebra $K[T]$. These identifications extend to a $K$-linear
equivalence $\Ev_K^{\Natural[1]^{op}}\simeq_K K[T]$-{\bf Mod}, where
$K[T]$-{\bf Mod} denotes the $K$-linear abelian category of
$K[T]$-modules.
}
\end{ex}

\paragraph{\S 4.2. Relation between both constructions.}
As mentioned before, even for finite categories $\Cc$, $\Ev_K^{\Cc^{op}}$ and
$\ev_K[\Cc]$ are generally non
equivalent. This is easily understood because $\Ev_K^{\Cc^{op}}$
is always abelian, while $\ev_K[\Cc]$ is not (see \S 3.6).

To make precise the relation between both constructions, observe that, among
the 
objects in $\Ev_K^{\Cc^{op}}$, we have the {\it representable} functors,
isomorphic to the functors $K[{\rm
  Hom}_{\Cc}(-,X)]:\Cc^{op}\To\Ev_K$ for
some object $X$ in $\Cc$ (note that, if the hom-sets of $\Cc$ are finite, such
functors actually take values in $\ev_K$). Then,
we have the following:

\begin{thm}
For any category $\Cc$ (resp. category $\Cc$ whose hom-sets are finite),
$\ev_k[\Cc]$ is $K$-linear equivalent to the $K$-linear additive subcategory of
$\Ev_K^{\Cc^{op}}$ (resp. of $\ev_K^{\Cc^{op}}$)  
generated by the representable functors (see Definition~\ref{definicio_base}).
\end{thm} 

\begin{proof}
For short, let $F_X$ stand for the functor $K[{\rm
  Hom}_{\Cc}(-,X)]$. Then, define a $K$-linear 
  functor $E:\ev_K[\Cc]\To\Ev_K^{\Cc^{op}}$ as
  follows:
\begin{itemize}
\item on objects: $E(X_1,\ldots,X_r)=F_{X_1}\oplus\cdots\oplus F_{X_r}$
  for   $r\geq 1$ and $E(\emptyset)=F_0$, the constant zero functor.
\item on morphisms: for any $f_{i'i}\in{\rm Hom}_{\Cc}(X_i,X'_{i'})$,
  let $A(f_{i'i}):(X_1,\ldots,X_r)\To(X'_1,\ldots,X'_{r'})$ 
  be the morphism with all entries equal to zero except the
  $(i',i)$-entry, which is equal to $f_{i'i}$. Then, define
  $E(A):F_{X_1}\oplus\cdots\oplus F_{X_r}\rightarrow 
F_{X'_1}\oplus\cdots\oplus F_{X'_{r'}}$ as the natural transformation whose
$Y$-component $E(A)_Y:F_{X_1}(Y)\oplus\cdots\oplus F_{X_r}(Y)\Rightarrow
F_{X'_1}(Y)\oplus\cdots\oplus F_{X'_{r'}}(Y)$ is the linear map described by
the $r'\times r$ matrix 
with entries $(E(A)_Y)_{j'j}:K[{\rm Hom}_{\Cc}(Y,X_j)]\To K[{\rm
  Hom}_{\Cc}(Y,X'_{j'})]$ given by
\begin{equation} \label{accio_E_sobre_morfismes}
(E(A)_Y)_{j'j}=\left\{ \begin{array}{ll} f_{i'i}\circ -, & \mbox{if $j=i$ and
      $j'=i'$} \\ 0, & \mbox{otherwise}
\end{array} \right.
\end{equation}
The morphisms $A(f_{i'i})$, for all $f_{i'i}\in{\rm
  Hom}_{\Cc}(X_i,X'_{i'})$ and all
$(i',i)\in\{1,\ldots,r'\}\times\{1,\ldots,r\}$, constitute a linear
basis of ${\rm
  Hom}_{\ev_K[\Cc]}((X_1,\ldots,X_r),(X'_{1},\ldots,X'_{r'}))$ and this action
of $E$ is extended $K$-linearly to arbitrary morphisms between
both sequences. 
\end{itemize}
It is easily checked that these assignments are
functorial. We only need to prove that it is a fully faithful functor.

Let us first see that the linear
map $E_{(X),(X')}:K[{\rm Hom}_{\Cc}(X,X')]\To{\rm Nat}(F_X,F_{X'})$
defined by (\ref{accio_E_sobre_morfismes}) is an isomorphism for any length
one sequences $(X),(X')$ of $\ev_K[\Cc]$. On the one hand,
we have a set bijection 
$$
{\sf Yon}:{\rm
  Hom}_{\Cc}(X,X')\To{\rm 
  Nat}({\rm Hom}_{\Cc}(-,X),{\rm Hom}_{\Cc}(-,X'))
$$
mapping $f:X\To X'$ to the
natural transformation $\sigma(f)$ with $Y$-component $\sigma(f)_Y=f\circ
-$ (Yoneda lemma). On the other hand, we have a linear map
$$
\Phi:K[{\rm
  Nat}({\rm Hom}_{\Cc}(-,X),{\rm Hom}_{\Cc}(-,X'))]\To{\rm
  Nat}(F_X,F_{X'})
$$
given by $\Phi(\sigma)=\sigma^K$, with $\sigma^K=1_{K[-]}\circ\sigma$, for all
$\sigma\in{\rm
  Nat}({\rm Hom}_{\Cc}(-,X),{\rm Hom}_{\Cc}(-,X'))$. The images
$\{\sigma^K\}_{\sigma}$ are linearly independent vectors of ${\rm 
  Nat}(F_X,F_{X'})$ (and hence, $\Phi$ is injective). Indeed, let
$\sigma_1\ldots,\sigma_n$ be arbitrary natural transformations from ${\rm
  Hom}_{\Cc}(-,X)$ to ${\rm Hom}_{\Cc}(-,X')$, with
$\sigma_i\neq\sigma_j$ if $i\neq j$. Note that $\sigma_i$ is
completely given by the morphism $(\sigma_i)_X({\id}_X)$ (Yoneda
lemma once more), so that the maps 
$(\sigma_1)_X({\id}_X),\ldots,(\sigma_n)_X({\id}_X):X\To X'$ are pairwise
different. Then, if $\sum_{i=1}^n\lambda_i\ \sigma_i^K:F_X\Rightarrow F_{X'}$
is the zero natural transformation, we have in particular that
$$
\sum_{i=1}^n\lambda_i\ (\sigma^K_i)_X({\id}_X)=\sum_{i=1}^n\lambda_i\
(\sigma_i)_X({\id}_X)=0, 
$$
in $K[{\rm Hom}_{\Cc}(X,X')]$. Hence, $\lambda_i=0$ for all
$i=1,\ldots,n$. Furthermore, given any natural transformation
$\tau:F_X\Rightarrow F_{X'}$, suppose that
$$
\tau_X({\id}_X)=\sum_{i=1}^n\lambda_i\ f_i
$$
where $f_i\in {\rm Hom}_{\Cc}(X,X')$, and define $\sigma_i:{\rm
  Hom}_{\Cc}(-,X)\Rightarrow {\rm Hom}_{\Cc}(-,X')$ 
by $(\sigma_i)_X({\id}_X)=f_i$ for all $i=1,\ldots,n$. Then, it is easily
checked that $\tau=\sum_{i=1}^n\lambda_i\ \sigma_i^K$, so that $\Phi$ is also
surjective. Therefore, $\Phi$ is an isomorphism of
vector spaces and it is immediately seen that the composite
\begin{align*}
K[{\rm Hom}_{\Cc}(X,X')]\stackrel{K[{\sf Yon}]}{\longrightarrow}& K[{\rm
  Nat}({\rm Hom}_{\Cc}(-,X),{\rm Hom}_{\Cc}(-,X'))]
\\ &\stackrel{\Phi}{\longrightarrow}{\rm
  Nat}(K[{\rm Hom}_{\Cc}(-,X)],K[{\rm Hom}_{\Cc}(-,X')])
\end{align*}
coincides with the linear map $E_{(X),(X')}$. 

More generally, for nonzero objects
$(X_1,\ldots,X_r)$ and $(X'_1,\ldots,X'_{r'})$ of arbitrary lengths, it is
$$
{\rm
  Hom}_{\ev_K[\Cc]}((X_1,\ldots,X_r),(X'_1,\ldots,X'_{r'}))\cong\prod_{(i,i')}
K[{\rm Hom}_{\Cc}(X_i,X'_{i'})]
$$
while
$$
{\rm Nat}(F_{X_1}\oplus\cdots\oplus F_{X_r},F_{X'_1}\oplus\cdots\oplus
F_{X'_{r'}})\cong \prod_{(i,i')}{\rm Nat}(F_{X_i},F_{X'_{i'}}),
$$
Under these identifications, it follows from the definition of $E$ that
$$
E_{(X_1,\ldots,X_r),(X'_1,\ldots,X'_{r'})}=\prod_{(i,i')}E_{(X_i),(X'_{i'})}
$$
Hence, the linear maps $E_{(X_1,\ldots,X_r),(X'_1,\ldots,X'_{r'})}$, for any
$r,r'\geq 1$, are also
isomorphisms and $E$ is indeed fully faithful.
\end{proof}

\begin{ex} \label{exemple_M=N_continuacio} {\rm
If $\Cc=\Natural[1]$, there is a unique representable functor up to
isomorphism. Namely, $F_{\ast}=K[{\rm Hom}_{\Natural[1]}(-,\ast)]$. Under the
identification $\Ev_K^{\Natural[1]^{op}}\simeq_K K[T]$-{\bf Mod} (see
Example~\ref{exemple_M=N}), this functor corresponds to $K[T]$ as module over
itself. Hence, $\ev_K[\Natural[1]]$ can be identified with the full
subcategory $K[T]$-{\bf Mod}$_f$ of $K[T]$-{\bf Mod} with objects the free
$K[T]$-modules. Note 
that this subcategory, and hence $\ev_K[\Natural[1]]$, is non Karoubian. Thus,
if $P$ is any projective non free $K[T]$-module and $F$ is 
the free $K[T]$-module of which $P$ is a direct summand, so that $F\cong
P\oplus M$ for some $K[T]$-submodule $M$ of $F$, the projection $p:F\To F$ of
$F$ onto $P$ is a non split idempotent in $K[T]$-{\bf Mod}$_f$ (an idempotent
$e:X\To X$ in an additive category splits if and only if ${\id}_X-e$ has
kernel, and ${\id}_X-p$ has no kernel in $K[T]$-{\bf Mod}$_f$).
}
\end{ex}

In some special cases, both categories
$\ev_K[\Cc]$ and $\Ev_K^{\Cc^{op}}$ may in fact be equivalent, mimicking the
situation for vector spaces. For instance, this is clearly the case if $\Cc$
is a finite discrete category and 
hence, for the Kapranov and Voevodsky 2-vector spaces. As shown by the
previous example, however, this is not true in general.

\section{General linear 2-groups $\G\Lb(\ev_K[\Cc])$}

Recall that for
any bicategory $\Cgg$ (non necessarily a strict one) and any 
object $X$ of $\Cgg$, the category ${\sf
  Equiv}_{\Cgg}(X)$ with objects the autoequivalences $f:X\To X$ and
with morphisms all 2-isomorphisms between these is a 2-group (see \S 2.6).

We are interested in the case $\Cgg=\Dev_K$. By analogy with the case
of vector spaces, let us denote by
$\G\Lb(\V)$ the 2-group ${\sf Equiv}_{\Dev_K}(\V)$ corresponding to a
generalized 2-vector space $\V$ and
call it the {\sf general linear 2-group} of $\V$. The purpose of this
section is to compute $\G\Lb(\V)$ (up to equivalence) for a special
type of generalized 2-vector spaces which include Kapranov and
Voevodsky 2-vector spaces.

\paragraph{\S 5.1. $\G\Lb(\ev_K[\Cc])$ versus ${\sf Equiv}_{\cat_K}(K[\Cc])$
  and ${\sf Equiv}_{\cat}(\Cc)$.}
Computing $\G\Lb(\V)$ for an 
arbitrary generalized 2-vector space seems to be difficult. There are,
however, general results relating 
$\G\Lb(\ev_K[\Cc])$ to the 2-groups ${\sf Equiv}_{\cat_K}(K[\Cc])$ and ${\sf
  Equiv}_{\cat}(\Cc)$ which we
want to discuss first, before considering any particular case.

Let $\V=\ev_K[\Cc]$, with $\Cc$
an arbitrary category, and let $H_{\Cc}:{\bf End}_{\cat_K}(K[\Cc])\To{\bf
  End}_{\Dev_K}(\ev_K[\Cc])$ be a functor mapping a $K$-linear functor
$\overline{F}:K[\Cc]\To K[\Cc]$ to some $K$-linear extension of the composite 
$K[\Cc]\stackrel{\overline{F}}{\To}
  K[\Cc]\stackrel{a_{K[\Cc]}}{\To}\ev_K[\Cc]$, 
  and a natural transformation $\overline{\tau}:\overline{F}\Rightarrow
  \overline{F}':K[\Cc]\To 
  K[\Cc]$ to the unique natural transformation $H_{\Cc}(\overline{\tau})$ such
  that $1_{a_{K[\Cc]}}\circ\overline{\tau}=H_{\Cc}(\overline{\tau})\circ
  1_{a_{K[\Cc]}}$ 
(cf. Proposition~\ref{propietat_universal_Add(C)}). There are various such
functors $H_{\Cc}$, but all of them are isomorphic because $K$-linear
extensions are unique up to isomorphism. They are clearly injective on
objects. In general, however, they are non
essentially surjective because a $K$-linear endomorphism
of $\ev_K[\Cc]$ can apply length one sequences to sequences of length
greater than one. For instance, if
$\Cc={\bf 1}$ (the terminal category), $K[{\bf 1}]$ is isomorphic to the one
object 
$K$-linear category $K[1]$ with $K$ as vector space of
endomorphisms, while $\ev_K[{\bf 1}]\simeq_K\ev_K$. Then,
${\bf End}_{\cat_K}(K[1])$ is a one object category, because a
$K$-linear 
functor $F:K[1]\To K[1]$ is nothing but a $K$-linear map
$f:K\To K$ and the condition of preservation of identities implies
that $f={\id}_K$
necessarily. In contrast, the set of isomorphism
classes of objects in ${\bf End}_{\Dev_K}(\ev_K)$ is in bijection with
the set $\Natural$ of natural numbers (see for ex. \cite{jE3}).

However, the following general result holds:

\begin{thm} \label{GL(C)_vs_Equiv(K[C])}
Let $\Cc$ be an arbitrary category. Then, any functor $H_{\Cc}:{\bf
End}_{\cat_K}(K[\Cc])\To{\bf End}_{\Dev_K}(\ev_K[\Cc])$ as above is a full
monoidal embedding. Moreover, if $\Cc$ is finite, it
restricts to an equivalence of monoidal categories (hence, an equivalence of
2-groups) 
$$
{\sf Equiv}_{\cat_K}(K[\Cc])\simeq \G\Lb(\ev_K[\Cc])
$$
\end{thm}
\begin{proof}
For any natural transformation
$\overline{\tau}_1,\overline{\tau}_2:\overline{F}\Rightarrow\overline{F}'$
between $K$-linear 
functors $\overline{F},\overline{F}':K[\Cc]\To K[\Cc]$, 
$H_{\Cc}(\overline{\tau}_1)=H_{\Cc}(\overline{\tau}_2)$ and the definition of
$H_{\Cc}$ on morphisms implies that
$1_{a_K[\Cc]}\circ\overline{\tau}_1=1_{a_K[\Cc]}\circ\overline{\tau}_2$ and
hence, $\overline{\tau}_1=\overline{\tau}_2$. Thus, $H_{\Cc}$ is always
faithful. It is also full because given 
$\sigma:H_{\Cc}(\overline{F})\Rightarrow H_{\Cc}(\overline{F}')$, the morphisms
$$ 
\overline{\tau}_X=\sigma_{(X)},\quad X\in{\rm Obj}(\Cc)
$$
define a natural transformation
$\overline{\tau}:\overline{F}\Rightarrow\overline{F}'$ 
such that $1_{a_{K[\Cc]}}\circ\overline{\tau}=\sigma\circ 1_{a_{K[\Cc]}}$ and
hence, we have $\sigma=H_{\Cc}(\overline{\tau})$. Furthermore, the functors
$H_{\Cc}(\overline{F}\overline{F}')$ and
$H_{\Cc}(\overline{F})H_{\Cc}(\overline{F}')$ both are $K$-linear
extensions of the functor 
$a_{K[\Cc]}\overline{F}\overline{F}':K[\Cc]\To\ev_K[\Cc]$, and
$H_{\Cc}({\id}_{K[\Cc]})$ and ${\id}_{\ev_K[\Cc]}$ both are $K$-linear 
extensions of $a_{K[\Cc]}:K[\Cc]\To\ev_K[\Cc]$. Consequently, there are
isomorphisms $H_{\Cc}(\overline{F}\overline{F}')\cong
H_{\Cc}(\overline{F})H_{\Cc}(\overline{F}')$ and
$H_{\Cc}({\id}_{K[\Cc]})\cong{\id}_{\ev_K[\Cc]}$. Let
$H_{\Cc,2}(\overline{F},\overline{F}'):
H_{\Cc}(\overline{F}\overline{F}')\Rightarrow
H_{\Cc}(\overline{F})H_{\Cc}(\overline{F}')$ and
$H_{\Cc,0}:H_{\Cc}({\id}_{K[\Cc]})\Rightarrow{\id}_{\ev_K[\Cc]}$ be the unique
natural isomorphisms such that
$1_{a_{K[\Cc]}\overline{F}\overline{F}'}=
H_{\Cc,2}(\overline{F},\overline{F}')\circ 1_{a_{K[\Cc]}}$ and $1_{a_{K[\Cc]}}=
H_{\Cc,0}\circ 1_{a_{K[\Cc]}}$
(cf. Proposition~\ref{propietat_universal_Add(C)}). Then, it is easy to check
that these
isomorphisms define a monoidal structure on $H_{\Cc}$. For instance, given any
$K$-linear endomorphisms $\overline{F},\overline{F}',\overline{F}''$ of
$K[\Cc]$, the following coherence condition needs to be checked:
$$
\left(H_{\Cc,2}(\overline{F},\overline{F}')\circ
  1_{H_{\Cc}(\overline{F}'')}\right)\cdot 
H_{\Cc,2}(\overline{F}\overline{F}',\overline{F}'')= 
\left(1_{H_{\Cc}(\overline{F})}\circ
  H_{\Cc,2}(\overline{F}',\overline{F}'')\right)\cdot 
H_{\Cc,2}(\overline{F},\overline{F}'\overline{F}'')
$$
Now, this equality holds if and only if the horizontal precomposites with
$1_{a_{K[\Cc]}}$ of both members are equal. But
\begin{align*}
\left[\left(H_{\Cc,2}(\overline{F},\overline{F}')\circ
  1_{H_{\Cc}(\overline{F}'')}\right)\cdot 
H_{\Cc,2}(\overline{F}\overline{F}',\overline{F}'')\right]\circ
  1_{a_{K[\Cc]}}&=\left(H_{\Cc,2}(\overline{F},\overline{F}')\circ
  1_{H_{\Cc}(\overline{F}'')a_{K[\Cc]}}\right)\cdot
  1_{a_{K[\Cc]}\overline{F}\overline{F}'\overline{F}''} 
\\ &=\left(H_{\Cc,2}(\overline{F},\overline{F}')\circ
  1_{a_{K[\Cc]}\overline{F}''}\right) \\
  &=1_{a_{K[\Cc]}\overline{F}\overline{F}'\overline{F}''} 
\end{align*}
and similarly
\begin{align*}
\left[\left(1_{H_{\Cc}(\overline{F})}\circ
    H_{\Cc,2}(\overline{F}',\overline{F}'')\right)\cdot  
H_{\Cc,2}(\overline{F},\overline{F}'\overline{F}'')\right]\circ
  1_{a_{K[\Cc]}}&=\left((1_{H_{\Cc}(\overline{F})}\circ
    H_{\Cc,2}(\overline{F}',\overline{F}'')\circ
  1_{a_{K[\Cc]}}\right)\cdot
  1_{a_{K[\Cc]}\overline{F}\overline{F}'\overline{F}''}  
\\ &=1_{H_{\Cc}(\overline{F})a_{K[\Cc]}\overline{F}'\overline{F}''} \\ 
  &=1_{a_{K[\Cc]}\overline{F}\overline{F}'\overline{F}''} 
\end{align*}
We leave to the reader checking the remaining coherence conditions and the
naturality of $H_{\Cc,2}(\overline{F},\overline{F}')$ in
$(\overline{F},\overline{F}')$. 

Suppose now that $\Cc$ is finite. Let $H_{\Cc}^0$ be the restriction of
$H_{\Cc}$ to ${\sf Equiv}_{\cat_K}(K[\Cc])$. It is a fully
faithful monoidal functor $H_{\Cc}^0:{\sf
  Equiv}_{\cat_K}(K[\Cc])\To\G\Lb(\ev_K[\Cc])$. For any $K$-linear functor 
$\breve{F}:\ev_K[\Cc]\To\ev_K[\Cc]$, let 
$F=\breve{F}\beta_{\Cc}:\Cc\To\ev_K[\Cc]$ be its restriction to
$\Cc$. $\breve{F}$ is obtained from $F$ by extending it with the 
help of some
biproduct functors and zero object in $\ev_K[\Cc]$.
Consequently, if $\breve{F}$ is an equivalence, the set of
image objects $\{F(X)\}_{X\in{\rm Obj}(\Cc)}$ generates $\ev_K[\Cc]$
additively. In particular, this set necessarily contains (up to isomorphism)
all indecomposable objects in $\ev_K[\Cc]$. But for a finite category $\Cc$,
all length one sequences 
are indecomposable (see Proposition~\ref{cas_C_finita}). It follows that
there exists a 
(unique) functor $\tilde{F}:\Cc\To K[\Cc]$, which uniquely extends to a
$K$-linear functor $\overline{F}:K[\Cc]\To K[\Cc]$, making commutative the
diagram 
$$
\xymatrix{
\ev_K[\Cc]\ar[rr]^{\breve{F}} & & \ev_K[\Cc] \\
& \Cc\ar[ru]^{F}\ar[ld]_{k_{\Cc}}\ar[rd]^{\exists_1\tilde{F}} & 
\\ K[\Cc]\ar[uu]^{a_{K[\Cc]}}\ar[rr]_{\exists_1\overline{F}} & &
K[\Cc]\ar[uu]_{a_{K[\Cc]}} 
}
$$
(in particular, $\overline{F}$ is an equivalence). By the
uniqueness up to isomorphism of the $K$-linear extensions, we conclude that
$\breve{F}\cong H_{\Cc}(\overline{F})$ and hence, $H^0_{\Cc}$ is also
essentially surjective. 
\end{proof}

Similarly, for any category $\Cc$, the functor ${\bf End}_{\cat}(\Cc)\To{\bf
End}_{\cat_K}(K[\Cc])$ mapping $F:\Cc\To\Cc$ to the unique
$K$-linear extension $\overline{F}$ of the composite
$\Cc\stackrel{F}{\To}\Cc\stackrel{k_{\Cc}}{\To}K[\Cc]$ and $\tau:F\Rightarrow
F'$ to the unique natural transformation 
$\overline{\tau}:\overline{F}\Rightarrow\overline{F}'$ such that
$\tau=\overline{\tau}\circ 1_{k_{\Cc}}$ is a monoidal non essentially
surjective embedding. But it is a non necessarily full functor now. Thus, in
the simplest case $\Cc={\bf 1}$, it is ${\bf
End}_{\cat}({\bf 1})\cong{\bf 1}$ while ${\bf
End}_{\cat_K}(K[{\bf 1}])\cong K[1]$ (any scalar $\lambda\in K$
defines a natural endomorphism of ${\id}_{K[{\bf 1}]}$). Furthermore, the
restriction ${\sf Equiv}_{\cat}(\Cc)\To{\sf Equiv}_{\cat_K}(K[\Cc])$ continues
to be neither full (${\sf Equiv}_{\cat_K}(K[{\bf
1}])\cong K^*[1]$ is not a terminal category) nor essentially surjective
(basically, because an arbitrary $K$-linear 
equivalence $K[\Cc]\To K[\Cc]$ need not map 
morphisms of $\Cc$ to morphisms also in $\Cc$).

\begin{ex} {\rm 
Let $\Cc$ be a group $G$ thought of as a category with only one object. Then,
if the restriction ${\sf Equiv}_{\cat}(\Cc)\To{\sf Equiv}_{\cat_K}(K[\Cc])$
really gives an essentially surjective functor, any
$K$-linear equivalence $K[\Cc]\To K[\Cc]$ should be isomorphic to the
$K$-linear extension of some equivalence $\Cc\To\Cc$. Now, a $K$-linear
equivalence $K[\Cc]\To K[\Cc]$ is nothing but a (unit preserving) algebra
automorphism of $K[G]$, while an equivalence $\Cc\To\Cc$ is just a
group automorphism of $G$. Furthermore, two algebra automorphisms
$\phi,\phi':K[G]\To K[G]$ define isomorphic $K$-linear equivalences if and only
if there exists a unit $u\in K[G]^*$ such that $\phi'(x)=u^{-1}\phi(x)u$ for
all $x\in K[G]$. In particular, if $G$ is abelian, they must be equal. But for
an abelian group $G$, an arbitrary automorphism of $K[G]$ does not restrict to
an automorphism of $G$. For ex., if $G=\Integer_2=\{\pm\}$, it is ${\rm
  Aut}_{{\bf Grp}}(\Integer_2)=1$, while ${\rm Aut}_{{\bf
    Alg}_K}(K[\Integer_2])\cong\Sigma_2$, the non trivial automorphism being
that which maps $(-)$ to $-(-)$ \footnote{The group ${\rm Aut}_{{\bf
    Alg}_K}(K[G])$ is computed in the next paragraph for an
arbitrary finite $G$ and an algebraically closed field $K$ whose
characteristic does not dived the order of $G$; see proof of 
Lemma~\ref{successio_escindida}.}.  }
\end{ex}

Therefore, the most general statement as regards 
the relation between $\G\Lb(\ev_K[\Cc])$ and ${\sf Equiv}_{\cat}(\Cc)$ reads as
follows: 

\begin{thm}
For any category $\Cc$, the composite functor ${\bf
  End}_{\cat}(\Cc)\hookrightarrow{\bf End}_{\cat_K}(K[\Cc])\To{\bf
  End}_{\Dev_K}(\ev_K[\Cc])$ restricts to a monoidal embedding ${\sf 
  Equiv}_{\cat}(\Cc)\hookrightarrow\G\Lb(\ev_K[\Cc])$. In particular, ${\sf
  Equiv}_{\cat}(\Cc)$ is equivalent to a (non full) sub-2-group of
  $\G\Lb(\ev_K[\Cc])$).  
\end{thm} 

For a finite category $\Cc$, this is to be thought of as an analog of the fact
that the group ${\rm Aut}(X)\cong\Sigma_n$ of automorphisms of a finite set
$X$ of cardinal 
$n$ is isomorphic to a subgroup (usually called the {\sf Weyl subgroup}) of
the general linear group ${\rm GL}(K[X])\cong{\rm GL}(n,K)$. This suggests
introducing the following 

\begin{defn}
For any finite category $\Cc$, the {\sf Weyl sub-2-group} of
$\G\Lb(\ev_K[\Cc])$ is the image of the previous monoidal embedding ${\sf
  Equiv}_{\cat}(\Cc)\hookrightarrow\G\Lb(\ev_K[\Cc])$.
\end{defn}  

It is not clear at all, however, that there exists some sort of analog for
$\G\Lb(\ev_K[\Cc])$ of the
Bruhat decomposition of the general linear groups ${\rm GL}(n,K)$.

\paragraph{\S 5.2. General linear 2-group of the generalized 2-vector space
  generated by a finite homogeneous groupoid.}
Recall that a groupoid $\Gg$ is a category equivalent to a
disjoint union of groups which are viewed as one object categories, i.e.,
$\Gg\simeq\sqcup_{i\in I}G_i[1]$ for some groups $G_i$. Let us call the
cardinal of $I$ the {\sf coarse size} of $\Gg$. We shall say that $\Gg$ is
{\sf homogeneous} when all groups $G_i$ are isomorphic to a given group
$G$, called the {\sf underlying group} of $\Gg$.

Suppose $\Gg$ is a finite homogeneous groupoid (i.e., finite coarse size and
finite underlying group). Examples include
all finite discrete categories 
$X[0]$ and all finite 2-groups $\G$, the first of coarse size equal to the
cardinal of $X$ and underlying group $G=1$, and the second of coarse size
equal to the cardinal of $\pi_0(\G)$ and underlying group $G=\pi_1(\G)$.

To simplify notation, we shall denote by $\G\Lb(\Gg)$ the general linear
2-group $\G\Lb(\ev_K[\Gg])$. The
purpose of this paragraph is to prove the following 

\begin{thm} \label{2-grup_lineal_general}
Let $K$ be an algebraically closed field and $\Gg$ a finite homogeneous
groupoid of coarse size $n$ and underlying group $G$. Suppose that the
order of $G$ is not divisible by the characteristic of $K$ (in particular,
this is the case if ${\sf char}(K)=0$). Then,
$\G\Lb(\Gg)$ is a split 2-group with 
\begin{align}
\pi_0(\G\Lb(\Gg))&
\cong\Sigma_n\times\left(\Sigma_{k_1}\times\cdots\times\Sigma_{k_s}\right)^n
\label{grup_pi_0}  \\
\pi_1(\G\Lb(\Gg))&\cong(K^*)^{rn} \label{grup_pi_1}
\end{align} 
where $\Sigma_p$ denotes the symmetric group on $p$ elements ($p\geq 1$), $r$
is the number of conjugacy classes of $G$ and $k_i\geq 1$, for all
$i=1,\ldots,s$, is the number of non equivalent
irreducible representations of $G$ of a given dimension $d_i$ (in particular,
$k_1+\cdots+k_s=r$). Furthermore, under these identifications, the action of  
$\pi_0(\G\Lb(\Gg))$ on $\pi_1(\G\Lb(\Gg))$ is given by
\begin{align} 
(\sigma,(\sigma_{1i})_{i=1}^s,\ldots,(\sigma_{ni})_{i=1}^s)\cdot
\left(\begin{array}{ccc} \lambda^{(1)}_{11} & \cdots &
    \lambda^{(1)}_{1n} \\ \vdots & & \vdots \\ \lambda^{(1)}_{k_11} & \cdots &
    \lambda^{(1)}_{k_1n} \\ \vdots & & \vdots \\  \lambda^{(s)}_{11} & \cdots &
    \lambda^{(s)}_{1n} \\ \vdots & & \vdots \\ \lambda^{(s)}_{k_s1} & \cdots &
    \lambda^{(s)}_{k_sn}  \end{array}\right)=
\left(\begin{array}{ccc} \lambda^{(1)}_{\sigma^{-1}_{11}(1)\sigma^{-1}(1)} &
    \cdots & 
    \lambda^{(1)}_{\sigma^{-1}_{n1}(1)\sigma^{-1}(n)} \\ \vdots & & \vdots \\
    \lambda^{(1)}_{\sigma^{-1}_{11}(k_1)\sigma^{-1}(1)} & \cdots & 
    \lambda^{(1)}_{\sigma^{-1}_{n1}(k_1)\sigma^{-1}(n)} \\ \vdots & & \vdots
    \\  \lambda^{(s)}_{\sigma^{-1}_{1s}(1)\sigma^{-1}(1)} & \cdots & 
    \lambda^{(s)}_{\sigma^{-1}_{ns}(1)\sigma^{-1}(n)} \\ \vdots & & \vdots \\
    \lambda^{(s)}_{\sigma^{-1}_{1s}(k_s)\sigma^{-1}(1)} & \cdots & 
    \lambda^{(s)}_{\sigma^{-1}_{ns}(k_s)\sigma^{-1}(n)}  \end{array}\right)
\nonumber \\ \mbox{} \nonumber \\ \mbox{} \label{accio_pi_0_sobre_pi_1}
\end{align}
for any $\sigma\in\Sigma_n$ and $\sigma_{qi}\in\Sigma_{k_i}$ for all
$i=1,\ldots,s$ and $q=1,\ldots,n$, and where we have identified the elements
$\Lambda\in(K^*)^{rn}$ with $r\times n$ 
matrices with entries $\lambda^{(i)}_{p_iq}\in K^*$. 
\end{thm}

Notice that for $G=1$ (hence, $r=1$) we indeed recover the general linear
2-groups of Kapranov 
and Voevodsky 2-vector spaces $\ev_K^n$, for which $\pi_0\cong\Sigma_n$,
$\pi_1\cong(K^*)^n$ and with $\Sigma_n$ acting on $(K^*)^n$ in the usual way,
i.e.
$$
\sigma\cdot(\lambda_1,\ldots,\lambda_n)=
(\lambda_{\sigma^{-1}(1)},\ldots,\lambda_{\sigma^{-1}(n)}) 
$$
(cf. \cite{jE4}, Proposition 6.3).

To prove the theorem, we shall first compute the homotopy groups $\pi_0$ and
$\pi_1$ 
of $\G\Lb(\Gg)$, next we shall determine the action of the first onto the
second and finally, we shall see that its classifying 3-cocycle $\alpha\in
Z^3(\pi_0(\G\Lb(\Gg)),\pi_1(\G\Lb(\Gg)))$ is cohomologically
trivial.

Recall that, for any $K$-algebra $A$, its group ${\rm Out}_{{\bf
Alg}_K}(A)$ of outer automorphisms is the quotient of the group 
${\rm Aut}_{{\bf Alg}_K}(A)$ of all its (unit preserving) algebra automorphisms
modulo the normal subgroup ${\rm Inn}_{{\bf Alg}_K}(A)$ of the inner ones,
i.e., of those of the form $\phi_u(x)=u^{-1}xu$ for some unit $u\in A^*$.

\begin{lem} \label{lema1_pi0}
For any finite homogeneous groupoid $\Gg$ of coarse size $n$ and underlying
group $G$, there is a group isomorphism
\begin{equation} \label{pi_0}
\pi_0(\G\Lb(\Gg))\cong \Sigma_n\times\left({\rm Out}_{{\bf
      Alg}_K}(K[G])\right)^n 
\end{equation}
In particular, $\pi_0(\G\Lb(\Gg))$ is a finite group.
\end{lem} 
\begin{proof}
Note first that for a groupoid $\Gg$ of the above kind it is
$K[\Gg]\simeq_K\sqcup_K^n K[G][1]$, where $\sqcup_K$ denotes the coproduct in
$\cat_K$ (see \S 2.4). Hence, a $K$-linear
equivalence of $K[\Gg]$ is completely determined by a permutation $\sigma\in
\Sigma_n$ giving the action on objects together with a collection of
$K$-algebra automorphisms $\phi_1,\ldots,\phi_n:K[G]\To K[G]$ giving the
action on the vector spaces of morphisms, and any such data
$(\sigma,\phi_1,\ldots,\phi_n)$ defines a $K$-linear equivalence
$F(\sigma,\phi_1,\ldots,\phi_n):K[\Gg]\To K[\Gg]$. Moreover, it is immediate to
check that the equivalences $F(\sigma,\phi_1,\ldots,\phi_n)$ and
$F(\sigma',\phi'_1,\ldots,\phi'_n)$ are isomorphic if and only if
$\sigma=\sigma'$ and there exists units $u_1,\ldots,u_n\in K[G]^*$ (the
components of an isomorphism) such that $\phi_i'(x)=u_i^{-1}\phi_i(x)u_i$ for
all $x\in K[G]$. The isomorphism (\ref{pi_0}) is then a consequence of
Theorem~\ref{GL(C)_vs_Equiv(K[C])}. As regards the last assertion, it follows
from a result due to Karpilovsky (see \cite{gK88}, Theorem
8.5.2), according to which the group ${\rm Out}_{{\bf
    Alg}_K}(K[G])$, for any $G$, is in bijection with the isomorphism classes
of $K[G\times G]$-modules whose underlying additive group is $K[G]$ and with
$K[G\times G]$ acting on it by
$$
\left(\sum_i\lambda_i(g_i,g'_i)\right)\
x=\sum_i\lambda_ig_ixf((g_i')^{-1}),\quad x\in K[G]
$$
for some $f\in{\rm Aut}_{{\bf Alg}_K}(K[G])$, and from the known fact (see
\cite{CR62}, Theorem 79.13) that for a finite group 
$G$, there are only finitely many isomorphism classes of such modules. 
\end{proof}

\begin{lem} \label{lema2_pi0}
Let $K$ be an algebraically closed field and $G$ a finite group whose order is
not divisible by the characteristic of $K$. Then, there is an isomorphism of
groups 
$$
{\rm Out}_{{\bf Alg}_K}(K[G])\cong\Sigma_{k_1}\times\cdots\times\Sigma_{k_s}
$$
where $k_i\geq 1$ is the number of non equivalent
irreducible representations of $G$ of a given dimension $d_i$, for all
$i=1,\ldots,s$ 
(in particular, $k_1+\cdots+k_s$ is the number of conjugacy classes of $G$).
\end{lem}
\begin{proof}
Under the assumptions on $K$ and on the order of $G$, it is well known (see
for ex. \cite{sW03}) that
there exists an algebra isomorphism $K[G]\cong M_{n_1}(K)\times\cdots\times
M_{n_r}(K)$,  
where $r$ is the number of conjugacy classes of $G$ and $n_1,\ldots,n_r$ are
the dimensions of the non equivalent irreducible representations of
$G$. Furthermore, it follows from Skolem-Noether theorem (see 
\cite{DK94}, Corollary 4.4.3) that all automorphisms of the algebra $M_n(K)$
are inner. Hence
\begin{align*}
{\rm Inn}_{{\bf Alg}_K}(M_{n_1}(K)\times\cdots\times M_{n_r}(K))&=
{\rm Inn}_{{\bf Alg}_K}(M_{n_1}(K))\times\cdots\times{\rm Inn}_{{\bf
    Alg}_K}(M_{n_r}(K)) \\
&= {\rm Aut}_{{\bf Alg}_K}(M_{n_1}(K))\times\cdots\times{\rm Aut}_{{\bf
    Alg}_K}(M_{n_r}(K))
\end{align*}
In general, however, the obvious embedding of ${\rm Aut}_{{\bf
Alg}_K}(M_{n_1}(K))\times\cdots\times{\rm Aut}_{{\bf Alg}_K}(M_{n_r}(K))$
into ${\rm Aut}_{{\bf Alg}_K}(M_{n_1}(K)\times\cdots\times M_{n_r}(K))$ is non
surjective, so that the quotient ${\rm Out}_{{\bf Alg}_K}(K[G])$ is non
trivial. To compute this quotient, let us denote by ${\bf I}_n$ the
identity $n\times 
n$ matrix and by {\bf 0} any zero matrix. Then, if $e_j=({\bf
  0},\ldots,\stackrel{j)}{{\bf 
  I}_{n_j}},\ldots,{\bf 0})$ ($j=1,\ldots,r$), the elements $e_1,\ldots,e_r$
are pairwise orthogonal central
idempotents of the product algebra $A=M_{n_1}(K)\times\cdots\times
M_{n_r}(K)$. Hence, any algebra 
automorphism $\phi:A\To A$ necessarily maps them to pairwise
orthogonal central idempotents of $A$. Since the center of $M_n(K)$ is
$Z(M_n(K))=K\ {\bf I}_n$, this means that 
$\phi(e_j)=\sum_{i=1}^r\lambda_{ij}e_i$ for some scalars
$\lambda_{ij}\in\{0,1\}$, $i,j=1,\ldots,r$ (note that $\phi(e_j)$ idempotent
implies that $\lambda_{ij}^2=\lambda_{ij}$). Moreover, since $\phi$
preserves the identity 
of $A$, we also have $\phi(e_1+\cdots+e_r)=e_1+\cdots+e_r$, from which it
follows that $\lambda_{ij}=\delta_{ij'}$, i.e., $\phi(e_j)=e_{j'}$ for some
$j'$ which depends on $j$. Together with the fact that, for any $N_j\in
M_{n_j}(K)$, it is
\begin{align*}
\phi({\bf 0},\ldots,N_j,\ldots,{\bf 0})&=\phi(({\bf 0},\ldots,N_j,\ldots,{\bf
  0})({\bf 0},\ldots,{\bf I}_{n_j},\ldots,{\bf 0}))
\\ &=\phi({\bf 0},\ldots,N_j,\ldots,{\bf
  0})\ \phi({\bf 0},\ldots,{\bf I}_{n_j},\ldots,{\bf 0})
\end{align*}
it follows that any automorphism $\phi$ of $A$
necessarily maps each factor $M_{n_j}(K)$ isomorphically onto some other
factor $M_{n_{j'}}(K)$. In particular, the subscript $j'$ for which
$\lambda_{ij}=\delta_{ij'}$ must be such that
$n_{j'}=n_j$. Inner or decomposable automorphisms $\phi\in{\rm 
Aut}_{{\bf Alg}_K}(M_{n_1}(K))\times\cdots\times{\rm Aut}_{{\bf 
Alg}_K}(M_{n_r}(K))$ correspond to the case $j'=j$ for all
$j=1,\ldots,r$. These will be the unique possible automorphisms of 
$A$ when the positive integers $n_1,\ldots,n_r$ are pairwise different. In
general, however, $G$ may have non equivalent irreducible representations of
the same dimension. Specifically, suppose we have $k_i$ non equivalent
irreducible representations of dimension $d_i$ for $i=1,\ldots,s$ (for
example, suppose that $n_1=\cdots=n_{k_1}=d_1$,
$n_{k_1+1}=\cdots=n_{k_1+k_2}=d_2$, etc.). In 
particular, we have $k_1+\cdots+k_s=r$. In this case, a generic automorphism of
$A=M_{d_1}(K)\times\stackrel{k_1)}{\cdots}\times 
M_{d_1}(K)\times\cdots\times M_{d_s}(K)\times\stackrel{k_s)}{\cdots}\times
M_{d_s}(K)$ will decompose in a unique way as the composite of a {\it
permutation} automorphism $\phi_{\sigma_1,\ldots,\sigma_s}$ given by
\begin{align} 
\phi_{\sigma_1,\ldots,\sigma_s}(N_{11},\ldots,N_{k_11},\ldots,&N_{1s},\ldots,N_{k_ss})= \nonumber
\\
&=(N_{\sigma_1(1)1},\ldots,N_{\sigma_1(k_1)1},\ldots,N_{\sigma_s(1)s},\ldots,N_{\sigma_s(k_s)s}) \label{automorfisme_permutacio} 
\end{align}
for some
$(\sigma_1,\ldots,\sigma_s)\in\Sigma_{k_1}\times\cdots\times\Sigma_{k_s}$,
followed by a decomposable automorphism $\phi_1\times\cdots\times\phi_r$. In
other words, by identifying $\Sigma_{k_1}\times\cdots\times\Sigma_{k_s}$ with
the above subgroup of permutation 
automorphisms of $A$, we conclude that ${\rm Aut}_{{\bf Alg}_K}(A)$ is the
semidirect product of the (normal) subgroup of inner automorphisms and
$\Sigma_{k_1}\times\cdots\times\Sigma_{k_s}$ and therefore, ${\rm Out}_{{\bf
    Alg}_K}(K[G])\cong\Sigma_{k_1}\times\cdots\times\Sigma_{k_s}$ as claimed. 
\end{proof}
Isomorphism (\ref{grup_pi_0}) of Theorem~\ref{2-grup_lineal_general} is now an
immediate consequence of Lemmas~\ref{lema1_pi0} and \ref{lema2_pi0}. Let us
now compute $\pi_1(\G\Lb(\Gg))$.

\begin{lem} \label{lema1_pi1}
For any finite homogeneous groupoid $\Gg$ of coarse size $n$ and underlying
group $G$, there is an isomorphism of abelian groups
\begin{equation} \label{pi_1}
\pi_1(\G\Lb(\Gg))\cong Z(K[G]^*)^n
\end{equation}
where $Z(K[G]^*)$ denotes the center of $K[G]^*$. 
\end{lem} 
\begin{proof}
By Theorem~\ref{GL(C)_vs_Equiv(K[C])}, $\pi_1(\G\Lb(\Gg))\cong\pi_1({\sf
  Equiv}_{\cat_K}(K[\Gg]))={\rm Aut}({\id}_{K[\Gg]})$. Now, a 
natural automorphism of ${\id}_{K[\Gg]}$ (actually, of any $F:K[\Gg]\To
  K[\Gg]$) is given by invertible elements
$u_1,\ldots,u_n$ (its components) in $K[G]$, and naturality further
requires that the $u_q$ belong to the center of
$K[G]^*$. Moreover, composition of automorphisms clearly corresponds to the
product in $Z(K[G]^*)^n$.
\end{proof}

\begin{lem} \label{lema2_pi1}
Let $K$ be an algebraically closed field and $G$ a finite group whose order is
not divisible by the characteristic of $K$. Then
$$
Z(K[G]^*)\cong (K^*)^r
$$
where $r$ is the number of conjugacy classes of $G$.
\end{lem}
\begin{proof}
The result is a direct consequence of the isomorphism of algebras $K[G]\cong
M_{n_1}(K)\times\cdots\times M_{n_r}(K)$ and the fact that
$Z(M_n(K))=K\ {\bf I}_n$. 
\end{proof}
Combining Lemmas~\ref{lema1_pi1} and \ref{lema2_pi1} we readily get isomorphism
(\ref{grup_pi_1}) of Theorem~\ref{2-grup_lineal_general}. Let us now prove
that, with the above identifications, the action is indeed given by
(\ref{accio_pi_0_sobre_pi_1}). 

\begin{lem}
For any finite homogeneous groupoid $\Gg$ of coarse size $n$ and underlying
group $G$, the action of $\pi_0(\G\Lb(\Gg))\cong\Sigma_n\times({\rm Out}_{{\bf
Alg}_K}(K[G]))^n$ on $\pi_1(\G\Lb(\Gg))\cong Z(K[G]^*)^n$
is given by 
\begin{equation} \label{accio}
(\sigma,[\phi_1],\ldots,[\phi_n])\cdot (u_1,\ldots,u_n)=
(\phi_{\sigma^{-1}(1)}(u_{\sigma^{-1}(1)}),\ldots,
\phi_{\sigma^{-1}(n)}(u_{\sigma^{-1}(n)})) 
\end{equation}
for any representatives $\phi_1,\ldots,\phi_n$ of
$[\phi_1],\ldots,[\phi_n]\in{\rm Out}_{{\bf Alg}_K}(K[G])$.
\end{lem}
\begin{proof}
Let us identify $\G\Lb(\Gg)$ with ${\sf Equiv}_{\cat_K}(K[\Gg])$. Then, by
definition, given $[F]\in\pi_0(\G\Lb(\Gg))$ and $\tau\in\pi_1(\G\Lb(\Gg))$, it
is  
$$
[F]\cdot \tau=\gamma_F^{-1}(\delta_F(\tau))
$$
for any representative $F:K[\Gg]\To K[\Gg]$ of $[F]$ (see \S 2.7). Now,
identifying ${\rm Aut}(F)$ with $Z(K[G]^*)^n$ as above, it is easy to see that
\begin{align*}
\delta_{F(\sigma,\phi_1,\ldots,\phi_n)}(u_1,\ldots,u_n)&=
(\phi_1(u_1),\ldots,\phi_n(u_n)) \\ 
\gamma_{F(\sigma,\phi_1,\ldots,\phi_n)}(u_1,\ldots,u_n)&=
(u_{\sigma(1)},\ldots,u_{\sigma(n)})
\end{align*}
from which (\ref{accio}) readily follows (note that the action so defined is
indeed independent of the representatives $\phi_i$ because the $u_i$ are
central). 
\end{proof}

Equation (\ref{accio_pi_0_sobre_pi_1}) follows now from (\ref{accio}) by
making the appropriate identifications. Thus, as discussed above, when $K[G]$
is identified 
with the corresponding product algebra $A=M_{n_1}(K)\times\cdots\times
M_{n_r}(K)$, each equivalence class
$[\phi_q]\in{\rm Out}_{{\bf Alg}_K}(A)$ in (\ref{accio}) can be identified
with an element $(\sigma_{q1},\ldots,\sigma_{qs})\in
\Sigma_{k_1}\times\cdots\times\Sigma_{k_s}$, a representative of $[\phi_q]$
being then the permutation automorphism
$\phi_{\sigma_{q1},\ldots,\sigma_{qs}}$ defined by
(\ref{automorfisme_permutacio}). At the same time, each $u_q\in Z(A^*)$,
$q=1,\ldots,n$, can 
be identified with an element $(\lambda^{(1)}_{1q},\ldots,\lambda^{(1)}_{k_1q},\ldots,\lambda^{(s)}_{1q},\ldots,\lambda^{(s)}_{k_sq})\in(K^*)^r$,
corresponding in fact to $(\lambda^{(1)}_{1q}{\bf
  I}_{d_1},\ldots,\lambda^{(1)}_{k_1q}{\bf 
  I}_{d_1},\ldots,\lambda^{(s)}_{1q}{\bf
  I}_{d_s},\ldots,\lambda^{(s)}_{k_sq}{\bf 
  I}_{d_s})\in Z(A^*)$. With these identifications, it is straightforward
checking that (\ref{accio}) translates into (\ref{accio_pi_0_sobre_pi_1}). 

Let us finally see that $\G\Lb(\Gg)$ is split. Note first the following
general result:

\begin{lem} \label{Equiv_vs_Aut}
For any finite category $\Cc$, there is an equivalence of 2-groups ${\sf
  Equiv}_{\cat_K}(K[\Cc])\simeq{\sf Aut}_{\cat_K}(K[\Cc])$ and hence, an
  equivalence of 2-groups $\G\Lb(\ev_K[\Cc])\simeq{\sf Aut}_{\cat_K}(K[\Cc])$.
\end{lem}
\begin{proof}
To prove the first assertion, it is enough to see that ${\sf
Equiv}_{\cat_K}(K[\Cc])={\sf Aut}_{\cat_K}(K[\Cc])$ for any finite {\it
skeletal} category $\Cc$. The claimed equivalence follows then from the fact
that any category is equivalent to a skeletal 
one. Let $\Cc$ be skeletal and let $E:K[\Cc]\To K[\Cc]$ be any
$K$-linear equivalence. In particular, there exists a $K$-linear functor
$\overline{E}:K[\Cc]\To K[\Cc]$ and a natural isomorphism
$\tau:\overline{E}E\Rightarrow {\id}_{K[\Cc]}$. Since $\Cc$ is finite, it
follows from Proposition~\ref{K-estabilitat} that $K[\Cc]$ is also
skeletal. Hence, 
$\overline{E}E(X)\cong X$ implies $\overline{E}E(X)=X$. Then, if 
$\widetilde{E}:K[\Cc]\To K[\Cc]$ is the $K$-linear functor uniquely defined by
$\widetilde{E}(X)=X$ for any object $X$ of $\Cc$ and $\widetilde{E}(f)=\tau_Y
f\tau_X^{-1}$ for any morphism $f:X\To Y$ in $\Cc$, it is easily checked that
$\widetilde{E}\overline{E}$ is a strict inverse of $E$. The last assertion
follows from Theorem~\ref{GL(C)_vs_Equiv(K[C])}.
\end{proof}
Let us now consider the case $\Cc$ is a finite homogeneous groupoid $\Gg$,
which we may assume it is skeletal. By
the previous Lemma, to prove that $\G\Lb(\Gg)$ is split it is enough to see
that ${\sf Aut}_{\cat_K}(K[\Gg])$ is split. But this is a strict 2-group and
hence, Proposition~\ref{2grup_escindit} and the subsequent remark can be
applied. Now, if $\Gg$ is skeletal, we have an strict equality $K[\Gg]=\sqcup^n
K[G][1]$ and hence
$$
|{\sf Aut}_{\cat_K}(\sqcup^n K[G][1])|=\Sigma_n\times\left({\rm
    Aut}_{{\bf Alg}_K}(K[G])\right)^n,
$$
while
$$
\pi_0({\sf Aut}_{\cat_K}(\sqcup^n K[G][1]))=\Sigma_n\times\left({\rm
    Out}_{{\bf Alg}_K}(K[G])\right)^n
$$
Therefore, the split character of ${\sf Aut}_{\cat_K}(K[\Gg])$ readily follows
from the next result:

\begin{lem} \label{successio_escindida}
Let $K$ be an algebraically closed field and $G$ a finite group whose order is
not divisible by the characteristic of $K$. Then, the exact sequence of groups
$1\To{\rm Inn}_{{\bf Alg}_K}(K[G])\To{\rm Aut}_{{\bf Alg}_K}(K[G])\To{\rm
  Out}_{{\bf Alg}_K}(K[G])\To 1$ splits.
\end{lem} 
\begin{proof}
It has already been shown that, under the above hypothesis on $K$ and $G$, the
group ${\rm Aut}_{{\bf Alg}_K}(K[G])$ is the semidirect product of ${\rm
Inn}_{{\bf Alg}_K}(K[G])$ and the subgroup of the so called permutation
automorphisms, which is isomorphic to a direct product of symmetric groups
(see proof of Lemma~\ref{lema2_pi0}). 
\end{proof}

\begin{cor}
For any finite 2-group $\G$, the general linear 2-group of the 2-vector space
it generates is split and with homotopy groups
\begin{align*}
\pi_0(\G\Lb(\G))&\cong\Sigma_n\times\Sigma_p^n \\
\pi_1(\G\Lb(\G))&\cong(K^*)^{pn}
\end{align*}
with $n$ and $p$ the cardinals of $\pi_0(\G)$ and $\pi_1(\G)$, respectively.
\end{cor}
\begin{proof}
$G$ is a finite homogeneous groupoid of coarse size $n$ and underlying group
$\pi_1(\G)$, which is abelian. 
\end{proof}

\section{Final comments}

The notion of generalized 2-vector space introduced in this work includes
Kapranov and Voevodsky 2-vector spaces, defined as a special kind of
$\ev_K$-module 
category (see \cite{KV94}). It is then worth comparing the notion of
generalized 2-vector space with the general notion of $\ev_K$-module category. 
In this sense, it is tedious but not difficult to see that any generalized
2-vector space $\ev_K[\Cc]$ has a ``canonical'' $\ev_K$-module category
structure, with $\ev_K$ acting on $\ev_K[\Cc]$ by
$$
V\odot S=\left\{ \begin{array}{cl} (S,\stackrel{n)}{\ldots},S) & \mbox{if ${\sf
          dim}V=n$} \\ \emptyset & \mbox{if ${\sf dim}V=0$} \end{array}\right.
$$
for any vector space $V$ and any object $S$ of $\ev_K[\Cc]$, and
$$
f\odot A=\left(\begin{array}{ccc} \alpha_{11}A & \cdots & \alpha_{1n}A
    \\ \vdots & & \vdots \\ \alpha_{n'1}A & \cdots & \alpha_{n'n}A
  \end{array} \right) 
$$
for any linear map $f:V\To V'$ and morphism $A:S\To S'$, where
$(\alpha_{i'i})$ is the matrix of $f$ in previously chosen linear bases of $V$
and $V'$ \footnote{This is a special case of the ``canonical'' $\ev_K$-module
  category 
  structure that can be defined on any $K$-linear additive category $\Aa$ once
  we choose particular biproduct functors of all orders and a zero object in
  $\Aa$ as well as a linear basis in each vector space.}. This can be seen as
the analog of the canonical $K$-linear structure on the sets $K[X]$. However,
it is unlikely that an arbitrary $\ev_K$-module
category is equivalent to a generalized 2-vector space
equipped with such a $\ev_K$-module category structure, as it happens with
vector spaces. Consequently, our notion of generalized 2-vector space should
still be 
thought of as a particular kind of $\ev_K$-module category, althought of a less
restrictive kind than Kapranov and Voevodsky 2-vector spaces.

To finish this paper, it seems worth pointing out the drawbakcs our notion of
generalized 2-vector space has with respect to Kapranov and Voevodsky 2-vector
spaces. In particular, let us mention some nice properties Kapranov and
Voevodsky 2-vector spaces have that are lost we moving to generalized 2-vector
spaces. 
One such lost property has already been mentioned before. Namely, an arbitrary
generalized 2-vector space is not always a Karoubian category (see \S
3.6). Another important drawback concerns the property of having a dual
object. Thus, an arbitrary generalized 2-vector space has no dual object, at 
least in the usual sense of the term, except when it is a Kapranov and
Voevodsky 2-vector. Indeed, if for a given $K$-linear
additive category $\Aa$ there exists a $K$-linear additive category $\Aa^*$
and $K$-linear functors ${\sf ev}:\Aa^*\Box\Aa\To\ev_K$ and ${\sf
  coev}:\ev_K\To\Aa\Box\Aa^*$ such that
$({\id}_{\Aa}\Box {\sf ev})({\sf 
  coev}\Box{\id}_{\Aa})\cong{\id}_{\Aa}$ and 
$({\sf ev}\Box{\id}_{\Aa^*})({\id}_{\Aa^*}\Box {\sf coev})\cong{\id}_{\Aa^*}$
\footnote{Here, $\Box$ denotes the tensor 
  product pseudofunctor for the 2-category $\adcat_K$ of $K$-linear additive
  categories; see \cite{jE6}.}, it may be shown that $\Aa$ 
is necessarily a Kapranov and Voevodsky 2-vector space (this
result seems to be due to P. Schauenburg; a proof can be found in Neuchl's
thesis \cite{mN97}). If 
$\Aa\simeq_K\ev_K^n$, such a dual object indeed exists and it is given
by the natural candidate, i.e., the hom-category ${\bf
  Hom}_{\dev_K}(\Aa,\ev_K)$, which is again a Kapranov and Voevodsky
2-vector space of rank $n$, in complete analogy with the situation for vector
spaces. Such a drawback may be serious when trying to define a
Frobenius structure on the 2-algebra $\ev_K[\G]$ generated by a 2-group
$\G$, because in the zero dimensional setting such a structure makes explicit
use of duals. But it is also possible that a definition of Frobenius structure
on a 2-algebra (more 
generally, on any pseudomonoid in a monoidal 2-category) may exists which
makes no use of duals (actually, such a definition where duals do not appear
already exists in the context of algebras or, more generally, monoids in a
monoidal category). 
Finally, another important drawback of generalized 2-vector spaces is that,
for arbitrary categories $\Cc$ and $\Dd$, the category of 
morphisms between the corresponding generalized 2-vector spaces ${\bf
  Hom}_{\Dev_K}(\ev_K[\Cc],\ev_K[\Dd])$ may no longer be a generalized
2-vector space. For
instance, if $\Cc=M[1]$, for some monoid $M$, and $\Dd={\bf 1}$,
this category of morphisms is nothing but the category of finite
dimensional linear representations of $M$, which is not a generalized 2-vector
space for an arbitrary $M$.

\vspace{0.5 truecm}
\noindent{{\bf Acknowledgments.}} I would like to acknoweledge B. Toen
for his many valuable comments, in particular, for pointing out to me the
relation discussed in Section 4 
between the constructions $\ev_K[\Cc]$ and $\Ev_K^{\Cc^{op}}$.

\bibliographystyle{plain}
\bibliography{2evb}

\end{document}